\input amssym
\input miniltx
\input graphicx.sty
\input pictex
\input color.sty

  %
  \font \bbfive = bbm5
  \font \bbseven = bbm7
  \font \bbten = bbm10

  \font \eightbf = cmbx8
  \font \eighti = cmmi8 \skewchar \eighti = '177
  \font \eightit = cmti8
  \font \eightrm = cmr8
  \font \eightsl = cmsl8
  \font \eightsy = cmsy8 \skewchar \eightsy = '60
  \font \eighttt = cmtt8 \hyphenchar \eighttt = -1

  \font \sixi = cmmi6 \skewchar \sixi = '177
  \font \sixrm = cmr6
  \font \sixsy = cmsy6 \skewchar \sixsy = '60
  \font \tensc = cmcsc10

  \scriptfont \bffam = \bbseven
  \scriptscriptfont \bffam = \bbfive
  \textfont \bffam = \bbten

  \newskip \ttglue

  \def \eightpoint {\def \rm {\fam 0 \eightrm }\relax
  \textfont 0= \eightrm
  \scriptfont 0 = \sixrm \scriptscriptfont 0 = \fiverm
  \textfont 1 = \eighti
  \scriptfont 1 = \sixi \scriptscriptfont 1 = \fivei
  \textfont 2 = \eightsy
  \scriptfont 2 = \sixsy \scriptscriptfont 2 = \fivesy
  \textfont 3 = \tenex
  \scriptfont 3 = \tenex \scriptscriptfont 3 = \tenex
  \def \it {\fam \itfam \eightit }\relax
  \textfont \itfam = \eightit
  \def \sl {\fam \slfam \eightsl }\relax
  \textfont \slfam = \eightsl
  \def \bf {\fam \bffam \eightbf }\relax
  \textfont \bffam = \bbseven
  \scriptfont \bffam = \bbfive
  \scriptscriptfont \bffam = \bbfive
  \def \tt {\fam \ttfam \eighttt }\relax
  \textfont \ttfam = \eighttt
  \tt \ttglue = .5em plus.25em minus.15em
  \normalbaselineskip = 9pt
  \def \MF {{\manual opqr}\-{\manual stuq}}\relax
  \let \sc = \sixrm
  \let \big = \eightbig
  \setbox \strutbox = \hbox {\vrule height7pt depth2pt width0pt}\relax
  \normalbaselines \rm }

  \def \setfont #1{\font \auxfont =#1 \auxfont}
  \def \withfont #1#2{{\setfont{#1}#2}}

  %

  \def \TRUE {Y}
  \def \FALSE {N}
  \def \EMPTY {}

  \def \ifundef #1{\expandafter \ifx \csname #1\endcsname \relax }

  \def \undefrule{\kern 2pt \vrule width 5pt height 5pt depth 0pt \kern 2pt}
  \def \UndefLabels{}
  \def \possundef #1{\ifundef {#1}\undefrule {\eighttt #1}\undefrule
    \global \edef \UndefLabels{\UndefLabels#1\par }
  \else \csname #1\endcsname \fi }

  %

  \newcount \secno \secno = 0
  \newcount \stno \stno = 0
  \newcount \eqcntr \eqcntr = 0

  \ifundef {showlabel} \global \def \showlabel {\FALSE} \fi  
  \ifundef {auxwrite} \global \def \auxwrite {\TRUE } \fi
  \ifundef {auxread} \global \def \auxread {\TRUE } \fi

  \def \define #1#2{\global \expandafter \edef \csname #1\endcsname {#2}}
  \long \def \error #1{\medskip \noindent {\bf ******* #1}}
  \def \fatal #1{\error{#1\par Exiting...}\end }

  \def \advseqnumbering {\global \advance \stno by 1 \global \eqcntr =0}

  \def \current {\ifnum \secno = 0 \number \stno \else \number \secno \ifnum \stno = 0 \else .\number \stno \fi \fi}

  \begingroup \catcode `\@=0 \catcode `\\=11 @global @def @textbackslash {\} @endgroup
  \def \space { }

  %
  \def \deflabel #1#2{%
    \if\TRUE\showlabel \hbox {\sixrm [[ #1 ]]} \fi
    \ifundef {#1PrimarilyDefined}%
      \define{#1}{#2}%
      \define{#1PrimarilyDefined}{#2}%
      \if\TRUE\auxwrite \immediate \write 1 {\textbackslash newlabel {#1}{#2}}\fi
    \else
      \edef \old {\csname #1\endcsname}%
      \edef \new {#2}%
      \if \old \new \else \fatal{Duplicate definition for label ``{\tt #1}'', already defined as ``{\tt \old}''.}\fi
      \fi}

  \def \label #1 {\deflabel {#1}{\current }}

  \def \equationmark #1 {\ifundef {InsideBlock}
	  \advseqnumbering
	  \eqno {(\current )}
	  \deflabel {#1}{\current }
	\else
	  \global \advance \eqcntr by 1
	  \edef \subeqmarkaux {\current .\number \eqcntr }
	  \eqno {(\subeqmarkaux )}
	  \deflabel {#1}{\subeqmarkaux }
	\fi }

  \def \split #1.#2.#3.#4;{\global \def \parone {#1}\global \def \partwo {#2}\global \def \parthree {#3}\global \def \parfour {#4}}
  \def \NA {NA}
  \def \ref #1{\split #1.NA.NA.NA;(\possundef {\parone }\ifx \partwo \NA \else .\partwo \fi )}

  %
  \newcount \bibno \bibno = 0

  \def \Bibitem #1 #2; #3; #4 \par{\smallbreak
    \global \advance \bibno by 1
    \item {[\possundef{#1}]} #2, {``#3''}, #4.\par
    \ifundef {#1PrimarilyDefined}\else
      \fatal{Duplicate definition for bibliography item ``{\tt #1}'', already defined in ``{\tt [\csname #1\endcsname]}''.}
      \fi
	\ifundef {#1}\else
	  \edef \prevNum{\csname #1\endcsname}
	  \ifnum \bibno=\prevNum \else
		\error{Mismatch bibliography item ``{\tt #1}'', defined earlier (in aux file ?) as ``{\tt \prevNum}'' but should be
	``{\tt \number\bibno}''.  Running again should fix this.}
		\fi
	  \fi
    \define{#1PrimarilyDefined}{#2}%
    \if\TRUE\auxwrite \immediate\write 1 {\textbackslash newbib {#1}{\number\bibno}}\fi}

  \def \jrn #1, #2 (#3), #4-#5;{{\sl #1}, {\bf #2} (#3), #4--#5}
  \def \Article #1 #2; #3; #4 \par{\Bibitem #1 #2; #3; \jrn #4; \par}

  \def \references {\begingroup \bigbreak \eightpoint \centerline {\tensc References} \nobreak \medskip \frenchspacing }

  %

  \catcode `\@=11
  \def \c@itrk #1{{\bf \possundef {#1}}} 
  \def \c@ite #1{{\rm [\c@itrk{#1}]}}
  \def \sc@ite [#1]#2{{\rm [\c@itrk{#2}\hskip 0.7pt:\hskip 2pt #1]}}
  \def \du@lcite {\if \pe@k [\expandafter \sc@ite \else \expandafter \c@ite \fi }
  \def \cite {\futurelet\pe@k \du@lcite }
  \catcode `\@=12

  %
  \def \Headlines #1#2{\nopagenumbers
    \headline {\ifnum \pageno = 1 \hfil
    \else \ifodd \pageno \tensc \hfil \lcase {#1} \hfil \folio
    \else \tensc \folio \hfil \lcase {#2} \hfil
    \fi \fi }}

  \def \title #1{\medskip\centerline {\withfont {cmbx12}{\ucase{#1}}}}

  \long \def \Quote #1\endQuote {\begingroup \leftskip 35pt \rightskip 35pt
\parindent 17pt \eightpoint #1\par \endgroup }
  \long \def \Abstract #1\endAbstract {\vskip 1cm \Quote \noindent #1\endQuote }

  \def \Note #1{\footnote {}{\eightpoint #1}}
  \def \Date #1 {\Note {\it Date: #1.}}

  \def\today{\ifcase\month\or January\or February\or March\or April\or May\or June\or July\or August\or September\or
    October\or November\or December\fi \space\number\day, \number\year}
  \def\hoje{\number\day/\number\month/\number\year}

  \newcount\auxone   \newcount\auxtwo   \newcount\auxthree
  \def\fulltime{\auxone=\time \auxtwo=\time \divide \auxone by 60 \auxthree=\auxone \multiply \auxthree by 60 \advance
    \auxtwo by -\auxthree \hoje\ \ \ifnum \auxone <10 0\fi\number\auxone :\ifnum \auxtwo<10 0\fi\number\auxtwo}

  \def \part #1#2{\vfill \eject \null \vskip 0.3\vsize
    \withfont{cmbx10 scaled 1440}{\centerline{PART #1} \vskip 1.5cm \centerline{#2}} \vfill\eject }

  %


  \def \fix {\smallskip \noindent $\blacktriangleright $\kern 12pt}
  \def \iskip {\medskip\noindent}

  \def \ucase #1{\edef \auxvar {\uppercase {#1}}\auxvar }
  \def \lcase #1{\edef \auxvar {\lowercase {#1}}\auxvar }

  \def \emph #1{{\it #1}}

  \def \section #1 \par {\global \advance \secno by 1 \stno = 0
    %
    \goodbreak \bigbreak
    \noindent {\bf \number \secno .\enspace #1.}
    \nobreak \medskip \noindent }

  \def \state #1 #2\par {\begingroup \def \InsideBlock {} \medbreak \noindent \advseqnumbering {\bf \current .\enspace
#1.\enspace \sl #2\par }\medbreak \endgroup }

  \def \definition #1\par {\state Definition \rm #1\par }


  \newcount \CloseProofFlag

  \long \def \Proof #1\endProof {\begingroup \def \InsideBlock {} \global \CloseProofFlag=0
     \medbreak \noindent {\it Proof.\enspace }#1
     \ifnum \CloseProofFlag=0 \hfill $\endproofmarker $ \looseness = -1 \fi
     \medbreak \endgroup}

  \def \quebra #1{#1 $$$$ #1}
  \def \explica #1#2{\mathrel {\buildrel \hbox {\sixrm #1} \over #2}}
  \def \explain #1#2{\explica{\ref{#1}}{#2}}  
  
  \def \=#1{\explain {#1}{=}}

  \def \pilar #1{\vrule height #1 width 0pt}

  \newcount \fnctr \fnctr = 0
  \def \fn #1{\global \advance \fnctr by 1
    \edef \footnumb {$^{\number \fnctr }$}%
    \footnote {\footnumb }{\eightpoint #1\par \vskip -10pt}}

  \def \text #1{\hbox {#1}}

  %
  \def \iItem {\smallskip }
  \def \Item #1{\smallskip \item {{\rm #1}}}
  \newcount \zitemno \zitemno = 0

  \def \izitem {\global \zitemno = 0}
  \def \iItemize {\izitem}
  \def \zitemplus {\global \advance \zitemno by 1 \relax }
  \def \rzitem {\romannumeral \zitemno }
  \def \rzitemplus {\zitemplus \rzitem } 
  \def \zitem {\Item {{\rm (\rzitemplus )}}}
  \def \iItem {\zitem }
  
  \def \zitemmark #1 {\deflabel {#1}{\current.\rzitem } \deflabel {Local#1}{\rzitem}}
  \def \iItemmark #1 {\zitemmark{#1} }

  \newcount \nitemno \nitemno = 0
  
  \def \nitem {\global \advance \nitemno by 1 \Item {{\rm (\number \nitemno )}}}

  \newcount \aitemno \aitemno = -1
  \def \boxlet #1{\hbox to 6.5pt{\hfill #1\hfill }}
  
  \def \aitemconv {\ifcase \aitemno a\or b\or c\or d\or e\or f\or g\or
h\or i\or j\or k\or l\or m\or n\or o\or p\or q\or r\or s\or t\or u\or
v\or w\or x\or y\or z\else zzz\fi }
  \def \aitem {\global \advance \aitemno by 1\Item {(\boxlet \aitemconv )}}
  \def \aitemmark #1 {\deflabel {#1}{\aitemconv }}

  \def \Case #1:{\medskip \noindent {\tensc Case #1:}}

  %
  \def \<{\left \langle \vrule width 0pt depth 0pt height 8pt }
  \def \>{\right \rangle }
  \def \({\big (}
  \def \){\big )}
  
  \def \and {\hbox {\quad and \quad }}

  \def \IFF {\kern 7pt\Leftrightarrow \kern 7pt}
  \def \IMPLY {\kern 7pt \Rightarrow \kern 7pt}
  \def   \for #1{\quad \forall \,#1}
  \def \endproofmarker {\square } 
  \def \"#1{{\it #1}\/} \def \umlaut #1{{\accent "7F #1}}
  \def \inv {^{-1}}
  \def \*{\otimes }
  \def \caldef #1{\global \expandafter \edef \csname #1\endcsname {{\cal #1}}}
  \def \mathcal #1{{\cal #1}}
  \def \bfdef #1{\global \expandafter \edef \csname #1\endcsname {{\bf #1}}}
  \bfdef N \bfdef Z \bfdef C \bfdef R
  
  \def \exists{\mathchar"0239\kern1pt }

  %

  \if \TRUE \auxread
    %
    \IfFileExists {\jobname.aux}{\input \jobname.aux}{\null} \fi
    \if \TRUE \auxwrite \immediate \openout 1 \jobname .aux \fi

  \def \close {\if \EMPTY \UndefLabels \else
      \message {*** There were undefined labels ***} \iskip
      ****************** \ Undefined Labels: \tt \par \UndefLabels
      \fi
    \if \TRUE \auxwrite \closeout 1 \fi
    \par \vfill \supereject \end }
  \def \bye {\close }
  %

  \def \inpatex #1.tex{\input #1.atex}

  \def \Caixa #1{\setbox 1=\hbox {$#1$\kern 1pt}\global \edef \tamcaixa {\the \wd 1}\box 1}
  \def \caixa #1{\hbox to \tamcaixa {$#1$\hfil }}


  \def\prod      {\mathop{\mathchoice{\hbox{$\mathchar "1351$}}{\mathchar "1351}{\mathchar "1351}{\mathchar "1351}}}
  \def\sum       {\mathop{\mathchoice{\hbox{$\mathchar "1350$}}{\mathchar "1350}{\mathchar "1350}{\mathchar "1350}}}
  \def\bigcup    {\mathop{\mathchoice{\raise 1pt \hbox{$\mathchar "1353$}}{\mathchar "1353}{\mathchar "1353}{\mathchar "1353}}}
  \def\bigcap    {\mathop{\mathchoice{\raise 1pt \hbox{$\mathchar "1354$}}{\mathchar "1354}{\mathchar "1354}{\mathchar "1354}}}

  \def\bigvee    {\mathop{\mathchoice{\hbox{$\mathchar "1357$}}{\mathchar "1357}{\mathchar "1357}{\mathchar "1357}}}


  \font \gothicfont =eufb10

  \def \Imply {\ \mathrel {\Rightarrow }\ }
  \def \Iff {\ \mathrel {\Leftrightarrow }\ }
  \def \|{\mathrel {|}}
  
  \def \arw #1{\ {\buildrel #1 \over \longrightarrow }\ }
  \def \tabrule {\noalign {\hrule }}

  \def \I {{\cal I}}
  \def \P {{\cal P}}
  
  \def \E {{\cal E}} 
  \def \S  {{\cal S}} 

  \def \goth #1{\hbox {\gothicfont #1}}
  \newdimen\largura \largura=5pt
  \def \gothEhat{\goth E\kern-\largura\widehat{\vrule height 6.5pt width \largura depth -10pt}}
  \def \hull {\goth H(S)}
  \def \ehull {\goth E(S)}

  \def \tight {{\hbox {\sixrm tight}}}
  \def \sixrmbox #1{{\hbox {\sixrm #1}}}
  \def \ess {\sharp}  

  \def\X{{\cal X}}  
  \def \subX {{\kern-1pt\scriptscriptstyle \X}}
  \def \SX {S_\subX}
  \def \tSX {\tilde \SX }
  \def \ehullX {\goth E(\SX )}
  \def \hullX {\goth H(\SX )}
  \def \specX {\gothEhat (\SX)}

  \def \eone {\goth E_1(S)}  

  \def\supop{^{_\sixrmbox{op}}}
  \def\subinf{_{^\infty}}

  \def\spec      {\gothEhat (S)}
  \def\varspec #1{\gothEhat #1(S)}
  \def\seone     {\varspec {_1}}
  \def\sinf      {\varspec {_{\scriptscriptstyle \infty}}}
  
  \def\smax      {\varspec {_\sixrmbox{max}}}
  \def\sess      {\varspec {_\sixrmbox{ess}}}
  \def\sinfop    {\varspec {\subinf\supop}}
  \def\sop       {\varspec {_\sixrmbox{op}}}

  \relax

  \def \id {\hbox {id}}
  \def \tS {\tilde S} 
  \def \sla {semilattice}
  
  \def \itmproof #1 {\medskip\noindent #1\enspace }
  \def \interior #1{\mathaccent'27#1}
  \definecolor {pink}{rgb}{1.0,0.21,0.7}
  \definecolor {grey}{rgb}{0.6,0.6,0.6}
  \definecolor {orange}{rgb}{1.0,0.5,0.0}
  \def \exists{\mathchar"0239\kern1pt }

\def \condef #1#2{\ifundef {#1}\expandafter \def \csname #1\endcsname{#2} \fi}

\def\writeIndex #1#2#3#4{\write #1 {\textbackslash #2 {#3}{#4}{\folio}}}

\def\printIndex#1#2#3#4{\noindent \hbox to #4{\hfill #1\hfill}\quad #2 \dotfill \quad #3 \pilar{11pt}\par}

\def\# #1; #2; {\writeIndex 3{symbol}{#1}{#2}}
\def \section #1 \par {\color{black}\global \advance \secno by 1 \stno = 0
  \goodbreak \bigbreak 
  \noindent {\bf \number \secno .\enspace #1.}  \nobreak \medskip \noindent
  \writeIndex 2{contents}{\number\secno.}{#1}}

\openout 3 symbols.aux

  \newcount \aitemno \aitemno = -1
  \def \boxlet #1{\hbox to 6.5pt{\hfill #1\hfill }}
  
  \def \aitemconv {\ifcase \aitemno a\or b\or c\or d\or e\or f\or g\or
h\or i\or j\or k\or l\or m\or n\or o\or p\or q\or r\or s\or t\or u\or
v\or w\or x\or y\or z\else zzz\fi }
  \def \aitem {\global \advance \aitemno by 1\Item {(\boxlet \aitemconv )}}
  \def \aitemmark #1 {\deflabel {#1}{\aitemconv }}


\centerline {\bf THE INVERSE HULL OF $0$\kern0.6pt-LEFT}
\smallskip
\centerline {\bf CANCELLATIVE SEMIGROUPS}
\medskip
\centerline {\tensc R. Exel and B. Steinberg}
\vskip 1cm


\Abstract Given a semigroup $S$ with zero, which is left-cancellative in the sense that $s t =s r\neq 0$ implies that $t
=r$, we construct an inverse semigroup called the inverse hull of $S$, denoted $\hull$.  When $S$ admits least common
multiples, in a precise sense defined below, we study the idempotent semilattice  of $\hull$, with a focus on
its spectrum.  When $S$ arises as the language semigroup for a subsift $X$ on a finite alphabet, we discuss the relationship between $\hull$ and
several C*-algebras associated to $X$ appearing in the literature.  \endAbstract

\openout 2 contents.aux

\section Introduction

The goal of this note is to announce a series of results about semigroups, together with applications to C*-algebras, whose proofs will appear in
later paper.  The theory of semigroup C*-algebras has a long history, beginning with Coburn's work \cite{CoOne} and \cite{CoTwo}
where the C*-algebra of the additive semigroup of the natural numbers is studied in connection to Toeplitz operators.
  In \cite{MurOne} G. Murphy generalized this construction to the positive cone of an ordered group, and later
to left-cancellative semigroups (\cite{MurTwo}, \cite{MurThree}).
  The C*-algebras studied by Murphy turned out to be too wild, even for nice looking semigroups such as ${\bf N}\times{\bf N}$, and this
prompted Li \cite{Li} to introduce an alternative C*-algebra for a left-cancellative semigroup.  By definition a
semigroup $S$ is said to be left-cancellative provided, for every $r,s,t\in S$, one has that
  $$
  s t =s r \Imply t =r.
  \equationmark Cancellative
  $$

Many interesting semigroups in the literature possess a zero element, namely an element $0$ such that
  $$
  s0=0s=0,
  $$
  for every $s$, and it is obvious that the presence of a zero prevents a semigroup from being left-cancellative.  In this
work we focus on $0$-left-cancellative semigroups, meaning that \ref{Cancellative} is required to hold only when the
terms in its antecedent are supposed to be nonzero.  This dramatically opens up the scope of applications including a
wealth of interesting semigroups, such as those arising from subshifts and, more generaly, languages over a fixed
alphabet.  This also allows for the inclusion of categories and  the semigroupoids of \cite{actions}, once the multiplication is
extended to all pairs of elements by setting undefined products to zero.

Starting with a $0$-left-cancellative semigroup $S$, the crucial point is to first build an inverse semigroup $\hull$,
which we call \emph{the inverse hull} of $S$, by analogy with \cite{CP},\cite{Cherubini}, from where one may invoke any of the now
standard constructions of C*-algebras from inverse semigroups, such as
  the tight C*-algebra \cite{actions} or
  Paterson's \cite{Paterson} universal C*-algebras.
  In fact this endeavor requires a lot more work regarding the passage from the original semigroup to its inverse hull,
rather than the much better understood passage from there to the C*-algebras.  Particularly demanding is the work geared
towards understanding the idempotent semilattice of $\hull$, which we denote by $\ehull$, as well as its spectrum.  By
a standard gadget $\ehull$ is put in correspondence with a subsemilattice of the power set of $S\setminus\{0\}$, whose
members we call the \emph{constructible sets}, by analogy with a similar concept relevant to Li's work in \cite{Li}.

Central to the study of the spectrum of $\ehull$ is the notion of \emph{strings},  which are motivated by the
description of the unit space of graph groupoids in terms of paths in the graph.

Regarding the problem of understanding the spectrum of $\ehull$, we believe the present work represents only a modest
beginning in a mammoth task lying ahead.  This impression comes from situations in which similar spectra have been more
or less understood, such as in \cite{infinoa} and in \cite{DokuchaExel}, illustrating the high degree of complexity
one should expect.

It is only in our final section that we return to considering C*-algebras where we discuss, from the present
perspective, how the Matsumoto and Carlsen-Matsumoto C*-algebras associated to a given subshift arise from the
consideration of the inverse hull of the associated language semigroup.  None of these correspond to the more well known
tight or Paterson's universal C*-algebras, but we show that they instead arise from reductions of the Paterson groupoid
to closed invariant subsets of its unit space which hitherto have not been identified.

\section Representations of semigroups

Let $S$ be a semigroup, namely a nonempty set equipped with an associative
operation.

A \emph {zero element}\/ for $S$ is a (necessarily unique) element $0\in S$,
satisfying
  $$
  s0 = 0s = 0, \for s\in S.
  $$

In what follows we will fix a semigroup $S$ possessing a zero element.

\definition \label DefRepre
  Let $\Omega $ be any set.  By a \emph {representation of $S$ on $\Omega $} we
shall mean any map
  $$
  \pi :S\to \I (\Omega ),
  $$
  \# $\textbackslash I (\Omega )$; Symmetric inverse  semigroup on the set $\Omega$;
  \# $\pi$; Representation of a semigroup;
  where $\I (\Omega )$ is the symmetric inverse
  semigroup\fn {The symmetric inverse semigroup on a set $\Omega $ is the
inverse semigroup formed by all partially defined bijections on $\Omega $.}  on
$\Omega $, such that
  \iItemize
  \iItem $\pi _0$ is the empty map on $\Omega $, and
  \iItem $\pi _s\circ \pi _t=\pi _{st}$, for all $s$ and $t$ in $S$.

Given a set $\Omega $, and any subset $X\subseteq \Omega $, let
  $\id _X$
  \# $\id _X$; Identity function on the set $X$;
  denote the identity function on $X$, so that $\id _X$ an element of $E\big (\I (\Omega )\big )$, the idempotent {\sla }
of $\I (\Omega )$.  One in fact has that
  $$
  E\big (\I (\Omega )\big ) = \{\id _X: X\subseteq \Omega \},
  $$
  \# $\textbackslash P(\Omega )$; The set of all subsets of $\Omega$;
  so we may identify $E\big (\I (\Omega )\big )$ with the meet {\sla } $\P
(\Omega )$ formed by all subsets of $\Omega $.

\definition \label IntroEsFs
  Given a representation $\pi $ of $S$, for every $s$ in $S$ we will denote the
domain of $\pi _s$ by $F^\pi _s$, and the range of $\pi _s$ by $E^\pi _s$, so
  \# $F^\pi _s$; Domain of $\pi_s$;
  \# $E^\pi _s$; Range of $\pi_s$;
  \# $F _s$; Simplified notation for domain of $\pi_s$;
  \# $E _s$; Simplified notation for range of $\pi_s$;
  that $\pi _s$ is a bijective mapping
  $$
  \pi _s:F^\pi _s \to E^\pi _s.
  $$
  When
  $$
  \Omega = \big (\bigcup _{s\in S}F^\pi _s \big ) \cup \big (\bigcup _{s\in S} E^\pi _s\big ),
  \equationmark EssSubset
  $$
  we will say that $\pi$ is an \emph {essential}\/ representation.
  We will moreover let
  \# $f^\pi _s$; Identity function on $F^\pi_s$;
  \# $e^\pi _s$; Identity function on $E^\pi_s$;
  $$
  f^\pi _s: = \pi _s\inv \pi _s = \id _{E^\pi _s}
  \and
  e^\pi _s:= \pi _s\pi _s\inv = \id _{F^\pi _s}.
  $$

Let us fix, for the time being, a representation $\pi $ of $S$ on $\Omega $.
Whenever there is only one representation in sight we will drop the superscripts
in $F^\pi _s$, $E^\pi _s$, $f^\pi _s$, and $e^\pi _s$, and adopt the simplified
notations $F_s$, $E_s$, $f_s$, and $e_s$.

The following may be proved easily.

\state Proposition \label Covar
  Given $s$ and $t$ in $S$, one has that
  \iItemize
  \iItem $\pi _se_t=e_{st}\pi _s$, and
  \iItem $f_t\pi _s=\pi _sf_{ts}$.

\definition \label GenSgANdConstr
  \iItemize
  \iItem The inverse subsemigroup of $\I (\Omega )$ generated by the set $\{\pi _s : s \in S\}$ will be denoted by
  $\I (\Omega ,\pi )$.
  \#   $\textbackslash I (\Omega ,\pi )$; The inverse semigroup  generated by the range of a representation of $\pi$;
  \iItem Given any $X\in \P (\Omega )$ such that $\id _X$ belongs to $E\big (\I (\Omega ,\pi )\big )$, we will say $X$
is a \emph {$\pi $-constructible subset}.
  \iItem The collection of all $\pi $-constructible subsets of $\Omega $ will be denoted by $\P (\Omega ,\pi )$.  In
symbols
  $$
  \P (\Omega ,\pi )= \big \{X\in \P (\Omega ) : \id _X \in E\big (\I (\Omega ,\pi )\big )\big \}.
  $$

  Observe that $E_s$ and $F_s$ are $\pi
$-constructible sets.  For the special case of $s=0$, we have $E_s=F_s=\emptyset$, so the empty set is $\pi $-constructible as well.

Since $\P (\Omega ,\pi )$ corresponds to the idempotent {\sla } of $\I
(\Omega ,\pi )$ by definition, it is clear that $\P (\Omega ,\pi )$ is a
{\sla }, and in particular the intersection of two $\pi $-constructible
sets is again $\pi $-constructible.

\section Cancellative semigroups

\definition
  Let $S$ be a semigroup containing a zero element.  We will say that $S$ is
\emph {$0$-left-cancellative} if,
for every $r,s,t\in S$,
  $$
  s t =s r \neq 0 \Imply t =r,
  $$
  and \emph {$0$-right-cancellative} if
  $$
  t s =r s \neq 0 \Imply t =r.
  $$
  If $S$ is both $0$-left-cancellative and $0$-right-cancellative, we will say
that $S$ is \emph {$0$-cancellative}.

\fix In what follows we will fix a $0$-left-cancellative semigroup $S$.  In a few
occasions we will also assume that $S$ is $0$-right-cancellative.

For any $s$ in $S$ we will let
  $$
  F_s = \{x\in S: sx\neq 0\},
  $$
  \#   $F_s$; Domains for the regular representation;
  and
  $$
  E_s = \{y\in S: y=sx\neq 0, \hbox { for some } x\in S\}.
  $$
  \#   $E_s$; Ranges for the regular representation;

Observe that the correspondence ``$x\to sx$'' gives a map from $F_s$ onto $E_s$,
which is one-to-one by virtue of $0$-left-cancellativity.

\definition \label DefineThetaS
  For every $s$ in $S$ we will denote by $\theta _s$ the bijective mapping given
by
  $$
  \theta _s: x\in F_s \mapsto sx\in E_s.
  $$
  \# $\theta$; Regular representation;

Observing that $0$ is neither in $F_s$, nor in $E_s$, we see that these are both
subsets of
  $$
  S':= S\setminus \{0\},
  \equationmark IntroSPrime
  $$
  \# $S'$; Set of all nonzero elements of $S$;
  so we may view $\theta _s$ as a partially defined bijection on $S'$, which is
to say that $\theta _s\in \I (S')$.  We should also notice that when $s=0$, both
$F_s$ and $E_s$ are empty, so $\theta _s$ is the empty map.

\state Proposition \label DefineRegRep
  The correspondence
  $$
  s \in S \mapsto \theta _s\in \I (S')
  $$
  is a representation of $S$ on $S'$, henceforth called the \emph {regular
representation of} $S$.

Regarding the notations introduced in \ref {IntroEsFs} in relation to the
regular representation, notice that
  $$
  F_s=F^\theta _s, \and E_s=E^\theta _s.
  $$

\state Definition \label RightReductive \rm
A semigroup $S$ is called \emph {right reductive} if it acts faithfully on the left of itself, that is, $sx=tx$ for all $x\in S$ implies $s=t$.

Of course every unital semigroup is right reductive.  If $S$ is a right reductive $0$-left cancellative semigroup, then
it embeds in $\I (S')$ via $s\mapsto \theta _s$.

Observe that if $S$ is $0$-right-cancellative, then a  single $x$ for which $sx=rx$, as long as this is nonzero, is enough to
imply that $s=t$.  So, in a sense, right reductivity is a weaker version of
$0$-right-cancellativity.

\definition
  The \emph {inverse hull} of $S$, henceforth denoted by $\hull $, is the
inverse subsemigroup of $\I (S')$ generated by the set $\{\theta _s : s \in
S\}$. Thus, in the terminology of \ref {GenSgANdConstr.i} we have
  \#   $\hull $; Inverse hull of $S$;
  $$
  \hull = \I (S',\theta ).
  $$

The reader should compare the above with the notion of \emph {inverse hull}
defined in \cite{CP}, \cite{Cherubini}.

The collection of $\theta $-constructible subsets of $S'$ is of special
importance to us, so we would like to give it a special notation:

\definition \label DefIdempSLA
  The idempotent {\sla } of $\hull $, which we will tacitly identify with
the {\sla } of $\theta $-constructible subsets of $S'$, will be denoted by
$\ehull $.
  Thus, in the terminology of \ref {GenSgANdConstr.iii} we have
  $$
  \ehull = \P (S',\theta ).
  $$
  \#   $\ehull $; Set of all $\theta$-constructible subsets of $S'$;

It will be important to identify some properties of $0$-left cancellative semigroups which will play a role later.

\state Proposition \label SomeCancelProps
Let $S$ be a $0$-left cancellative semigroup.
    \iItemize
    \iItem If $e$ is an idempotent element in $S$  and $s\in S\setminus \{0\}$, then $es\neq
0$ if and only if $es=s$, that is, $s\in eS\setminus \{0\}$.
    \iItem If $s\in S\setminus \{0\}$, then $s\in sS$ if and only if $se=s$ for a necessarily unique idempotent $e$.
    \iItem If $sS=S$ and $S$ is right reductive, then $S$ is unital and $s$ is invertible.

A semigroup $S$ is said to have \emph {right local units} if $S=SE(S)$, that is, for all $s\in S$, there exists an
idempotent element $e$ in $S$  with $se=s$.  A unital semigroup has right local units for trivial reasons.  If $S$ has right local units, then $sS=0$ implies that $s=0$.  From Proposition~\ref {SomeCancelProps} we obtain the following:

\state Corollary \label LocalUnits
Let $S$ be a $0$-left cancellative semigroup.  Then $S$ has right local units if and only if $s\in sS$ for all $s\in S$.

In a right reductive $0$-left cancellative semigroup, the idempotents are orthogonal to each other:

\state Proposition \label OrthogIdems
Let $S$ be right reductive $0$-left cancellative semigroup and suppose that $e\neq f$ are distinct idempotents of $S$.  Then $ef=0$.

If $S$ is a $0$-left cancellative, right reductive semigroup with right local units, then for $s\in S\setminus \{0\}$, we denote by
  $s^+$
  \#  $s^+$; Right local unit for $s$;
  the unique idempotent with $ss^+=s$.  If $S$ is unital, then $s^+=1$.  If $C$ is a left cancellative category, we can associate a semigroup $S(C)$ by letting $S(C)$ consist of the arrows of $C$ together with a
zero element $0$.  Products that are undefined in $C$ are made zero in $S(C)$. It is straightforward to check that
$S(C)$ is $0$-left cancellative, right reducitve and has right local units. If $f\colon c\to d$ is an arrow of $C$, then
$f^+=1_c$.

\section Examples

\label ExamplesSection
  If $G$ is any group, and $P$ is a subsemigroup of $G$, let $S=P\cup\{0\}$,   where $0$ is any element not belonging to $P$,
with the obvious multiplication operation.  Then $S$
 is clearly a $0$-cancellative semigroup.

A more elaborate example is as follows:
  let $\Lambda $ be any finite or infinite set, henceforth called the \emph{alphabet}, and let $\Lambda ^+$ be the free semigroup generated by $\Lambda $,
namely the set of all finite words in $\Lambda $ of positive length (and hence
excluding the empty word), equipped with the multiplication operation given by
concatenation.

Let $L$ be a \emph {language} on $\Lambda $, namely any nonempty subset of
$\Lambda ^+$.  We will furthermore assume that $L$ is \emph {closed under
prefixes and suffixes} in the sense that, for every $\alpha $ and $\beta $ in $\Lambda
^+$, one has
  $$
  \alpha \beta \in L \Imply \alpha \in L, \hbox { and } \beta \in L .
  $$

  Define a multiplication operation on
  $$
  S:=L\cup \,\{0\},
  \equationmark LanguageSGroup
  $$
  where $0$ is any element not belonging to $L$, by
  $$
  \def \quad { }
  \alpha \cdot \beta =\left \{\matrix {
  \alpha \beta , & \hbox { if $\alpha ,\beta \neq 0$, and $\alpha \beta \in L$},
\cr \pilar {12pt}
  0, & \hbox { otherwise.}\hfill }\right .
  $$
One may then prove that $S$ is a $0$-cancellative semigroup.

One important special case of the above example is based on \emph {subshifts}.
Given an alphabet $\Lambda $, as above, let $X\subseteq \Lambda ^{\bf N}$ be any
nonempty subset invariant under the left shift, namely the mapping $\sigma
:\Lambda ^{\bf N}\to \Lambda ^{\bf N}$ given by
  $$
  \sigma (x_1x_2x_3\ldots )=x_2x_3x_4\ldots
  $$

  Let $L\subseteq \Lambda ^+$ be the \emph {language of} $X$, namely the set of
all finite words occuring in some infinite word belonging to $X$. Then $L$ is
clearly closed under prefixes and suffixes, and hence we are back in the
conditions of the above example.

The fact that $X$ is invariant under the left shift is indeed superfluous, as
any nonempty subset $X\subseteq \Lambda ^{\bf N}$ would lead to the same
conclusion. However, languages arising from subshifts have been intensively
studied in the literature, hence the motivation for the above example. The semigroup associated to the language of a shift
was first studied in the early days of symbolic dynamics by Morse and Hedlund~\cite{MH}.

\definition (\cite {Munn})
  We will say that a semigroup $S$ with zero is \emph {categorical at zero} if,
for every $r,s,t\in S$, one has that
  $$
  rs\neq 0, \hbox { and } st\neq 0 \Imply rst\neq 0.
  $$

The semigroup associated to a category is categorical at zero, whence the name.
We may use the above ideas to produce a semigroup which is not categorical at
zero: take any nonempty alphabet $\Lambda $, let $L$ be the language consisting
of all words of length at most two, and let $S=L\cup \{0\}$, as described above.
If $a$, $b$ and $c$, are members of $\Lambda $, we have that $abc=0$, but $ab$
and $bc$ are nonzero, so $S$ is not categorical at zero.

In contrast, notice that if $S$ is the semigroup built as above from a Markov
subshift, then $S$ is easily seen to be categorical at zero.

\def\auto{\sigma}
Another interesting example is obtained from self-similar graphs \cite{ExelPardo}.
Let
    $G$ be a discrete group,
    $E=(E^0, E^1)$ be a graph with no sources,
    $\auto $ be an action of $G$ by graph automorphisms on $E$,
    and
  $$
  \varphi :G\times E^1 \to G
  $$
  be a one-cocycle for the restriction of $\auto $ to $E^1$, which moreover satisfies
  $$
  \auto _{\varphi (g,e )}(x ) = \auto _g(x )
  \for g \in G \for e \in E^1 \for x \in E^0.
  $$

Let $E^*$ be the set of all finite
  paths\fn{We adopt the functorial point of view so a path is a sequence $e_1\ldots e_n$ of edges, such that
$s(e_i)=r(e_{i+1})$, for all $i$.}
  on $E$, including the vertices, which are seen as paths of length zero, and put
  $$
  S=(E^*\times G)\cup\{0\}.
  $$
  Given nonzero elements $(\alpha,g)$ and $(\beta,h)$ in $S$, define
  $$
  (\alpha,g) (\beta,h) = \left\{ \matrix{ \big(\alpha \auto_g(\beta),  \varphi(g, \beta)h\big), &\hbox{if }s(\alpha)=r\big(\auto_g(\beta)\big),  \cr\cr
  0, & \hbox{otherwise.}\hfill} \right.
  $$

One may then show that $S$ is a categorical at zero, $0$-left-cancellative semigroup, which is $0$-right-cancella\-tive
if and only if $(G, E)$ is pseudo-free in the sense of \cite[Proposition 5.6]{ExelPardo}.

\section Least common multiples

We now wish to introduce a special class of semigroups, but  for this we must first consider the question of divisibility.

\definition
  Given $s$ and $t$ in a semigroup $S$, we will say that $s$
  \emph {divides}
  $t$, in symbols
  $$
  s\| t,
  $$
  \#  $s\| t$; $s$ divides $t$;
  or that $t$ is a \emph {multiple} of $s$, when either $s=t$, or there is some $u$ in $S$ such that
  $su=t$.  In other words,
  $$
  t\in\{s\}\cup sS.
  $$

We observe that division is a reflexive and transitive relation, so it may be seen as a (not necessarily anti-symmetric)
order relation.

For the strict purpose of simplifying the description of  the division relation, regardless of whether or not $S$ is
unital,  we shall sometimes employ the unitized semigroup
  $$
  \tS := S\cup\{1\},
  $$
  \#  $\textbackslash tS$; Unitized semigroup $S\cup\{1\}$;
  where $1$ is any element not belonging to $S$, made to act like a unit for $S$.  For every $s$ and $t$  in $S$ we therefore
have that
  $$
  s\|t \iff \exists u\in \tS , \ su=t.
  \equationmark UnitDivision
  $$

Having enlarged our  semigroup, we might as well extend the notion of divisibility:

\definition Given $v$ and $w$ in $\tS $,  we will say that $v\|w$ when there exists some $u$ in $\tS $, such that $vu=w$.

Notice that if $v$ and $w$ are in $S$, then the above notion of divisibility coincides with the previous one by
\ref {UnitDivision}.  Analysing the new cases where this extended divisibility may or may not apply, notice that:
  $$
  \matrix {
  \forall w\in\tS , &
  1\|w, \hfill \cr \pilar {20pt}
  \forall v\in\tS , &
  v\|1 \iff v=1.}
  \equationmark NewDivision
  $$

  The introduction of $\tS $ brings with it several pitfalls, not least because $\tS $ might not be
$0$-left-cancellative: when $S$ already has a unit, say $1_S$, then in the identity ``$s1_S=s1$'', we are not allowed to
left-cancel $s$, since $1_S\neq1$.
One should therefore be very careful when working with $\tS $.

\definition \label DefLCM
  Let $S$ be a semigroup and let $s,t\in S$.  We will say that an element $r\in S$ is a \emph {least common
multiple} for $s$ and $t$ when
  \iItemize
  \iItem $sS\cap tS = rS$,
  \iItem both $s$ and $t$ divide $r$.

Observe that when $S$ has right local units then $r\in rS$, by \ref {LocalUnits}, and hence condition \ref {DefLCM.i}
trivially implies \ref {DefLCM.ii}, so the former condition alone suffices to define least common multiples.
  However in a unitless semigroup condition  \ref {DefLCM.i} may hold while \ref {DefLCM.ii} fails.
Nevertheless, when $sS\cap tS=\{0\}$, then $0$ is a least common multiple for $s$ and $t$ because $s$ and $t$ always
divide
$0$.

Regardless of the existence of right local units, notice that condition \ref {DefLCM.ii} holds if and only if $r\tS \subseteq s\tS \cap t\tS $, and
therefore one
has that $r$ is a least common multiple for $s$ and $t$ if and only if
  $$
  sS\cap tS = rS \subseteq r\tS \subseteq s\tS \cap t\tS .
  \equationmark sStStilde
  $$

  One may of course think of alernative definitions for the concept of least common multiples, fiddling with the above
ideas in many different ways.  However \ref{DefLCM}
seems to be the correct choice, at least from the point of view of the theory we are about to develop.

\definition \label DefLCMSgrp
  We shall say that a semigroup $S$ \emph {admits least common multiples} if there exists a least common multiple for
each pair of elements of $S$.

The language semigroup of  \ref{LanguageSGroup} is easily seen to be an example of a
semigroup admitting least common multiples.

Another interesting example is obtained from the
quasi-lattice ordered groups of \cite {Nica}, which we would now like to briefly describe.

Given a group $G$ and a unital subsemigroup $P\subseteq G$, one defines a
partial order on $G$ via
  $$
  x\leq y \iff x\inv y \in P.
  $$
  The quasi-lattice condition says that, whenever elements $x$ and $y$ in $G$
admit a common upper bound, namely an element $z$ in $G$ such that $z\geq x$ and
$z\geq y$, then there exists a least common upper bound, usually denoted $x\vee
y$.

Under this situation, consider the semigroup $S = P\cup \{0\}$, obtained by
adjoining a zero to $P$.  Then, for every nonzero $s$ in $S$, i.e.~for $s$ in
$P$, one has that
  $$
  sS=\{x\in P: x\geq s\}\cup \{0\},
  $$
  so that the multiples of $s$ are precisely the upper bounds of $s$ in $P$,
including zero.

  If $t$ is another nonzero element in $S$, one therefore has that
  $s$ and $t$ admit a nonzero commun multiple if
and only if $s$ and $t$ admit a common upper bound in $P$, in which case $s\vee
t$ is a least common multiple of $s$ and $t$.

On the other hand, when $s$ and $t$ admit no common upper bound, then obviously
$s\vee t$ does not exist, but still $s$ and $t$ admit a least common multiple in
$S$, namely $0$.

Summarizing we have the following:

\state Proposition
  Let $(G,P)$ be a quasi-lattice ordered group.  Then the semigroup $S:=P\cup
\{0\}$ admits least common multiples.

\fix From this point on we will fix a $0$-left-cancellative, semigroup $S$
admitting least common multiples.

\state Proposition \label ExistLCMinSTilde
  Given $u$ and $v$ in $\tS $, there exists $w$ in $\tS $ such that
  \iItemize
  \iItem $uS\cap vS = wS$,
  \iItem both $u$ and $v$ divide $w$.

\Proof
  When $u$ and $v$ lie in $S$, it is enough to take $w$ to be a  (usual) least common multiple of $u$ and $v$.  On the other
hand, if $u=1$, one takes $w=v$, and if $v=1$, one takes $w=u$.
\endProof

Based on the above we may extend the notion of least common multiples to $\tS $, as follows:

\definition \label NewDefLCM
  Given $u$ and $v$ in $\tS $, we will say that an element $w$ in $\tS $ is a least common multiple of $u$ and $v$,
provided \ref {ExistLCMinSTilde.i-ii} hold.  In the exceptional case that $u=v=1$, only $w=1$ will be considered to be a
least common multiple of $u$ and $v$, even though there might be another $w$ in $S$ satisfying
\ref {ExistLCMinSTilde.i-ii}.

It is perhaps interesting to describe the exceptional situation above, where we are arbitrarily prohibiting by hand that an element of $S$  be
considered as a least common multiple of $1$ and itself, even though it would otherwise satisfy all of the required properties.
  If $w\in S$ is such an element, then
  $$
  wS = 1S\cap 1S=S,
  $$
  so, in case we also assume that $S$ is right-reductive, we deduce from \ref {SomeCancelProps.iii} that $S$ is unital
and $w$ is invertible.  Thus, in hindsight it might not have been such a good idea to add an external  unit to $S$ after
all.

On the other hand, when $s$ and $t$ lie in $S$, it is not hard to see that any least common multiple of $s$ and $t$ in the new sense of
\ref {NewDefLCM} must belong to $S$, and hence
it must also be a least common multiple in the old sense of \ref {DefLCM}.

Given a representation $\pi$ of $S$ on a set $\Omega$, we will now concentrate our  attention on giving
a concrete description
for the inverse semigroup $\I (\Omega ,\pi )$ defined in \ref {GenSgANdConstr},
provided $\pi $ satisfies certain special properties,
which we will now describe.

Initially   notice that if $s\|r$, then the range of $\pi_r$ is contained in the range of $\pi_s$ because
either $r=s$, or $r=su$, for some $u$ in $S$, in which case
  $
  \pi _r = \pi _s\pi _u.
  $
  \ So, using the notation introduced in \ref {IntroEsFs},
  $$
  E^\pi _r \subseteq E^\pi _s.
  $$
When $r$ is a common multiple of $s$ and $t$, it then follows that
  $$
  E^\pi _r \subseteq E^\pi _s \cap E^\pi _t.
  $$

\definition \label DefRespect
  A representation $\pi $ of $S$ is said to \emph {respect least common
multiples} if, whenever $r$ is a least common multiple of elements $s$ and $t$
in $S$, one has that
  $
  E^\pi _r = E^\pi _s \cap E^\pi _t.
  $

As an example, notice that the regular representation of $S$, defined in \ref
{DefineRegRep}, satisfies the above condition since the fact that $rS=sS\cap tS$
implies that
  $$
  E^\theta _r =
  rS\setminus \{0\} =
  (sS \cap tS)\setminus \{0\} =
  (sS\setminus \{0\}) \cap (tS\setminus \{0\}) =
  E^\theta _s \cap E^\theta _t.
  \equationmark RanLCM
  $$

\fix From now on we will fix a representation $\pi $ of $S$ on a set $\Omega $, assumed to respect least common multiples.
Since $\pi $ will be the only representation considered for a while, we will
  use the simplified notations
  $F_s$, $E_s$, $f_s$, and $e_s$.

There is a cannonical way to extend $\pi$ to $\tS $ by
setting
  $$
  F_1=E_1=\Omega, \and \pi_1=\id _\Omega.
  $$
  It is evident that $\pi$ remains a multiplicative map after  this extension.
Whenever we find it convenient we will therefore think of $\pi$ as defined on $\tS $ as above.   We will accordingly extend
the notations
  $f_s$ and $e_s$ to allow for any $s$ in $\tS $, in the obvious way.

\state Proposition \label ExtensionRespLCM
  Let $\pi$ be a representation of\/ $S$ on a set\/ $\Omega$.  If $\pi$ respects least common multiples then so does
its natural extension to $\tS $.  Precisely, if $u$ and $v$ are elements of $\tS $, and if $w\in\tS $ is a least common multiple
of $u$ and $v$, then
  $
  E _w = E _u \cap E _v.
  $

\definition
  Given a representation $\pi $ of $S$, and given any nonempty finite subset
$\Lambda \subseteq \tS $, we will let
  $$
  F^\pi _\Lambda = \bigcap _{u\in \Lambda } F^\pi _u, \and
  f^\pi _\Lambda = \prod _{u\in \Lambda }f^\pi _u.
  $$
  \#   $F^\pi _\Lambda $; Intersection of $F^\pi _u$, for $u$ in $\Lambda$;
  \#   $f^\pi _\Lambda $; Product of $f^\pi _u$, for $u$ in $\Lambda$;

When there is only one representation of $S$ in sight, as in the present moment, we will drop the superscripts and use the simplified notations
  $F_\Lambda $ and
  $f_\Lambda $.

We should remark that, since each $f _s$ is the identity map on $F _s$,
one has that $f _\Lambda $ is the identity map on $F _\Lambda $.

Also notice that, since $f_1=\id _\Omega$,
the presence of $1$ in $\Lambda$ has no effect in the sense that $f_\Lambda = f_{\Lambda\cup\{1\}}$, for every $\Lambda$.  Thus, whenever convenient
we may assume that $1\in\Lambda$.

As already indicated we are interested in obtaining a description for the inverse semigroup $\I (\Omega,\pi)$.  In that respect
it is interesting to observe that
most elements of the form
  $f_\Lambda $  belong to $\I (\Omega, \pi)$, but there is one exception,
namely when $\Lambda=\{1\}$.  In this case we have
  $$
  f_{\{1\}} = \id _\Omega,
  $$
  which may or may not lie in $\I (\Omega, \pi)$.  However, when $\Lambda\cap S\neq\emptyset $, then surely
  $$
  f_\Lambda \in  \I (\Omega, \pi).
  \equationmark fLambdaInHull
  $$

We should furthermore  remark that, whenever we are looking at a term of the form
  $\pi _uf_\Lambda \pi _v\inv $, we may assume that $u,v\in \Lambda $, because
  $$
  \pi _uf_\Lambda \pi _v\inv =
  \pi _u\pi _u\inv \pi _uf_\Lambda \pi _v\inv \pi _v\pi _v\inv \quebra =
  \pi _uf_u f_\Lambda f_v \pi _v\inv =
  \pi _uf_{\{u\}\cup \Lambda \cup \{v\}} \pi _v\inv ,
  \equationmark LambdaHasu
  $$
  so $\Lambda $ may be replaced by $\{u\}\cup \Lambda \cup \{v\}$
without altering the above term.
Moreover, as in \ref {fLambdaInHull}, observe that  when $\Lambda\cap S\neq\emptyset $, then
  $$
  \pi _uf_\Lambda \pi _v\inv \in  \I (\Omega, \pi).
  $$


\state Theorem \label GenFormHull
  Let $S$ be a $0$-left-cancellative semigroup admitting least common multiples.
  Also let $\pi $ be a representation of\/ $S$ on a set $\Omega $, assumed to
respect least common multiples.
  Then the generated inverse semigroup $\I (\Omega ,\pi )$ is given by
  $$
  \I (\Omega ,\pi ) =
  \big \{\pi _uf_\Lambda \pi _v\inv :
  \Lambda \subseteq \tS \hbox { is finite, } \Lambda\cap S\neq\emptyset , \hbox { and } u,v\in \Lambda
  \big \}.
  $$

With this we may describe the constructible sets in a very concrete way.

\state Proposition
  Under the assumptions of \ref {GenFormHull}, the $\pi $-constructible subsets
of\/ $\Omega $ are precisely the sets of the form
  $$
  X= \pi _u(F_\Lambda ) ,
  $$
  where $\Lambda \subseteq \tS $ is a finite subset, $\Lambda\cap S\neq\emptyset $,
and $u\in \Lambda $.

Recalling that the regular representation of $S$ respects least common multiples,
our last two results apply to give:

\state Corollary \label FormOfHull
  Let $S$ be a $0$-left-cancellative semigroup admitting least common multiples.
Then
  $$
  \hull =
  \big \{\theta _uf_\Lambda \theta _v\inv :
  \Lambda \subseteq \tS \hbox { is finite, } \Lambda\cap S\neq\emptyset , \hbox { and } u,v\in \Lambda
  \big \},
  $$
  and
  $$
  \ehull =
  \big \{uF_{\Lambda } ,\
  \Lambda \subseteq \tS \hbox { is finite, } \Lambda\cap S\neq\emptyset , \hbox { and } u\in \Lambda
  \big \}.
  $$

If $S$ is categorical at zero and has right local units, then a more precise normal form can be obtained.

\section Strings

This section is intended to introduce a device which will be highly useful in the study of the spectrum of $\ehull$.

\fix Throughout this section $S$ will be a fixed $0$-left-cancellative semigroup.

\definition \label DefineStrin
  A nonempty subset
  $\sigma\subseteq S$
  \# $\sigma$; A string;
  is said to be a \emph {string} in $S$, if
  \iItemize
  \iItem $0\notin\sigma$,
  \iItem for every $s$ and $t$ in $S$, if $s\|t$, and $t\in\sigma$, then $s\in\sigma$, \iItemmark StringHereditary
  \iItem for every $s_1$ and $s_2$ in $\sigma$, there is some $s$ in $\sigma$ such that $s_1\| s$, and $s_2\| s$.  \iItemmark StringDirected

\bigskip An elementary example of a string is the set of divisors of any
nonzero element $s$ in $S$,  namely,
  $$
  \delta_s= \{t\in S: t\| s\}.
  \equationmark StringDivisors
  $$
  \# $\delta_s$; The string of divisors of $s$;

Strings often contain many elements, but there are some exceptional strings consisting of a single semigroup element. To
better study these it is useful to introduce some terminology.

\definition \label DefinePrime
  Given $s$ in $S$ we will say that $s$ is:
  \iItemize
  \iItem \emph {prime}, if the only divisor of $s$ is $s$, itself, or, equivalently, if $\delta_s=\{s\}$,
  \iItem \emph {irreducible}, if there are no two elements $x$ and $y$ in $S$ such that $s=xy$, or,  equivalently,  if $s\notin S^2$.

It is evident that any irreducible element is prime, but there might be prime elements which are not irreducible.  For
example, in the semigroup $S=\{0,  s,  e\}$, with multiplication table given by

  \bigskip \hfill \vbox {\begingroup
  \def \tabrule {\noalign {\hrule }}
  \def \vr {\vrule height 12pt}
  \def \bx #1{\hbox to 20pt {\ \hfill #1\hfill }}
  \offinterlineskip
  \halign {
    \strut
    \vr #&\vr #&\vr #&\vr #\vr \cr
    \tabrule
   \bx {$\times$} & \bx {$0$} & \bx {$e$} & \bx {$s$}\cr \tabrule
   \bx {$0$} & \bx {$0$} & \bx {$0$} & \bx {$0$}\cr \tabrule
   \bx {$e$} & \bx {$0$} & \bx {$e$} & \bx {$0$}\cr \tabrule
   \bx {$s$} & \bx {$0$} & \bx {$s$} & \bx {$0$}\cr \tabrule
  }
  \endgroup }\hfill \null

  \bigskip \noindent one has that $s$ is prime but not irreducible because $s=se\in S^2$.

\state Proposition \label PrimeStrings A singleton $\{s\}$ is a string if and only if $s$ is prime.

\Proof If $s$ is prime then the singleton $\{s\}$ coincides with $\delta_s$, and hence it is a string.  Conversely, supposing that
$\{s\}$ is a string, we have by \ref {StringHereditary} that $\delta_s\subseteq\{s\}$, from where it follows that $s$ is prime.
  \endProof

\definition
  The set of all strings in $S$ will be denoted by $S^\star $.
  \# $S^\star $;   The set of all strings in $S$;

From now on our goal will be to define an action of $S$ on $S^\star $.

\state Proposition \label MappingStrings
  Let $\sigma$ be a string in $S$,  and let $r\in S$.  Then
  \iItemize
  \iItem if\/ $0$ is not in $r\sigma$, one has that
  \setbox 1\hbox {$\{t\in S: t\| rs, \hbox { for some } s\in\sigma\}$}
  \setbox 0\hbox {$r\inv *\sigma$}
  $$
  \hbox to \wd 0{\hfill $r*\sigma$}:= \copy 1
  $$
  \# $r*\sigma$; Product of the semigroup element $r$ by the string $\sigma$;
  is a string whose intersection with $rS$ is nonempty.
  \iItem If $\sigma$   is a string whose intersection with $rS$ is nonempty, then
  $$
  \copy 0:= \hbox to \wd 1 {$\{t\in S: rt\in\sigma\}$\hfill }
  $$
  \# $r\inv *\sigma$; Inverse product of the semigroup element $r$ by the string $\sigma$;
  is a string, and $0$ is not in $r(r\inv *\sigma)$.

It should be noted that, under the assumptions of \ref {MappingStrings.i}, one
has that
  $$
  r\sigma\subseteq r*\sigma,
  \equationmark ProdInForw
  $$
  and in fact $r*\sigma$ is the hereditary closure of $r\sigma$ relative to the order relation given by division.

We may then define a representation of $S$ on the set $S^\star $ of all strings in $S$,   as follows:

\state Proposition \label IntroStarAction
  For each $r$ in $S$, put
  \# $\theta^\star $; Representation of $S$ on $S^\star $;
  \# $F^\star _r$; Domain of $\theta^\star _r$;
  \# $E^\star _r$; Range of $\theta^\star _r$;
  $$
  F^\star _r=\{\sigma \in S^\star : r\sigma \not \ni 0\}, \and
  E^\star _r=\{\sigma\in S^\star : \sigma \cap rS \neq\emptyset \}.
  $$
  Also let
  $$
  \theta ^\star _r:F^\star _r\to E^\star _r
  $$
  be defined by $\theta ^\star _r(\sigma ) = r*\sigma $, for every $\sigma \in
F^\star _r$.  Then:
  \iItemize
  \iItem $\theta ^\star _r$ is bijective, and its inverse is the mapping defined
by
  $$
  \sigma \in E^\star _r\mapsto r\inv *\sigma \in F^\star _r.
  $$
  \iItem Viewing $\theta ^\star $ as a map from $S$ to $\I (S^\star )$, one has
that $\theta ^\star $ is a representation of $S$ on $S^\star $.

Useful alternative characterizations of $F^\star _r$ and $E^\star _r$ are as follows:

\state Proposition \label FRFstarR
  Given $r$ in $S$, and given any string $\sigma$ in $S^\star $, one has that:
  \iItemize \def \iff {\ \mathrel {\Leftrightarrow }\ }
  \iItem $\sigma\in F^\star _r \Iff \sigma\subseteq F^\theta_r$,
  \iItem $\sigma\in E^\star _r \Iff \sigma\cap E^\theta_r\neq\emptyset ,$
  \iItem $\sigma\in E^\star _r \Imply r\in\sigma$.  In addition, the converse holds provided $r\in rS$ (e.g.~if $S$  has right local units).

Recall from \ref {FormOfHull} that, when $S$ has least common multiples, every $\theta$-constructible subset of $S'$ has the
form $\theta_u(F^\theta_\Lambda)$, where $\Lambda\subseteq\tS $ is finite, $\Lambda\cap S\neq\emptyset $, and $u\in\Lambda$.
  By analogy this suggests that it might also be useful to have a characterization of $\theta^\star _u(F^\star _\Lambda)$ along the lines of
\ref {FRFstarR}.

\state Proposition \label FRFstarRTwo
  Let $\Lambda$ be a finite subset of $\tS $ having a nonempty intersection with $S$, and let $u\in\Lambda$.  Then
  $\theta^\star _u(F^\star _\Lambda)$
  consists precisely of the strings $\sigma$ such that
  $$
  \emptyset \neq \sigma\cap E^\theta_u \subseteq \theta_u(F^\theta_\Lambda).
  $$

After \ref {IntroStarAction} we now have two natural representations of $S$, namely the regular representation $\theta$
acting on $S'$, and $\theta^\star $ acting on $S^\star $.

\state Proposition The map
  $$
  \delta: s\in S'\mapsto\delta_s\in S^\star ,
  $$
  where $\delta_s$ is defined in \ref {StringDivisors}, is covariant relative to $\theta$ and
$\theta^\star $.

Observe that the union of an increasing family of strings is a string, so any
string is contained in a maximal one by Zorn's Lemma.

\definition
  The subset of $S^\star $ formed by all maximal strings will be denoted by $S^\infty$.
  \# $S^\infty$;    The set of all maximal strings;

  Our next result says that $S^\infty$ is invariant under $\theta^\star $.

\state Proposition \label ForwardInvariance
  For every $r$ in $S$,  and for every maximal string $\sigma$ in
  $F^\star _r$,
  one has that $\theta^\star _r(\sigma)$ is maximal.

Observe that the above result says that $S^\infty$ is invariant under each
$\theta^\star _r $, but not necessarily under ${\theta^\star _r}\inv $.

An example  to show that $S^\infty$ may indeed be non invariant under
${\theta^\star _r}\inv $ is as follows.  Consider the language $L$ on the alphabet $\Sigma=\{a, b\}$ given by
  $$
  L=\{a, b, aa, ba\}.
  $$
  Then,  $\sigma=\{b, ba\}$ is a maximal string, while ${\theta^\star _b}\inv =\{a\}$ is not maximal.

\definition By a representation of a given inverse semigroup $\S $ on a set $\Omega$ we shall mean any map
  $$
  \rho:\S \to\I (\Omega),
  $$
  such that $\rho(0)$ is the empty map, and for every $s$ and $t$ in $\S $, one has that $\rho(st)=\rho(s)\rho(t)$, and
$\rho(s\inv )=\rho(s)\inv $.

\state Proposition \label MapForISG
  Let $S$ be a $0$-left-cancellative semigroup admitting least common multiples.
  Then there exists a unique representation $\rho$ of $\hull $ on  $S^\star $,
  such that the following diagram commutes. \smallskip
  \hfill \beginpicture \setcoordinatesystem units <0.003truecm, -0.00250truecm>
  \put {$S$} at 0000 0000
  \put {$\I (S^\star )$} at 1100 0000
  \put {$\hull $} at 500 580
  \arrow <0.11cm> [0.3,1.2] from 150 000 to 850 000
  \put {$\theta^\star $} at 500 -130 \arrow <0.11cm> [0.3,1.2] from 100 100 to 400 450
  \put {$\theta$} at 150 350
  \arrow <0.11cm> [0.3,1.2] from 600 450 to 900 100
  \put {$\rho$} at 830 350
  \endpicture \hfill \null

Observing that a homomorphism of inverse semigroups must restrict to the
corresponding idempotent {\sla }s, we obtain the following:

\state Corollary \label IntroEpsilon
  Let $S$ be a $0$-left-cancellative semigroup admitting least common multiples.
  Then there exists a {\sla } representation
  $$
  \varepsilon:\ehull \to\P (S^\star ),
  $$
  \#   $\varepsilon$; Representation of $\ehull $ on $S^\star $;
  such that
  $$
  \varepsilon\big (\theta_u(F^\theta_\Lambda)\big )=\theta^\star _u(F^\star _\Lambda),
  $$
  whenever $\Lambda$ is a finite subset of $\tS $ intersecting $S$ nontrivially, and $u\in\Lambda$.

Observing that $E^\theta_r=\theta_r(F^\theta_r)$,  notice that
  $$
  \varepsilon(E^\theta_r) = \theta^\star _r(F^\star _r)=E^\star _r.
  $$

\section The spectrum of the semilattice of constructible sets

Let us now fix a $0$-left-cancellative semigroup $S$ admitting least common multiples.  It is our goal in this
section to study the spectrum of $\ehull $.

Recall that if $\E$ is a semilattice with zero, the spectrum of $\E$ is the set of all semilattice homomorphisms $\varphi:\E\to\{0,
1\}$, such that $\varphi(0)=0$.  Here $\{0, 1\}$ is equipped with its standard {\sla } structure $0<1$.

Considering the representation
  $$
  \varepsilon:\ehull \to\P (S^\star ),
  $$
  introduced in \ref {IntroEpsilon},
  and given $\sigma\in S^\star $, set
  $$
  \varphi_\sigma : X\in\ehull \mapsto [\sigma\in\varepsilon(X)] \ \in\ \{0, 1\}.
  $$
  \#  $\varphi_\sigma$; Character associated with the string $\sigma$;
  It is clear that $\varphi_\sigma$ is a {\sla } homomorphism, so
it is a character
  as long as it is nonzero.

The question of whether or not $\varphi_\sigma$ is nonzero evidently boils down to the existence of some constructible set $X$ for
which $\sigma\in\varepsilon(X)$.   By \ref{FormOfHull} it is easy to see that for every constructible set $X$,  there is some $r$ in $S$
such that $X\subseteq E_r$, or $X\subseteq F_r$.
Therefore $\varphi_\sigma=0$ if and only if $\sigma$ is never in any $E^\star_r$ nor in any $F^\star_r$,  that is,
  $$
  \sigma\not\subseteq F^\theta_r, \and \sigma\cap E^\theta_r=\emptyset,
  $$
  for every $r$ in $S$, by \ref{FRFstarR}.

The second condition above implies that every element in $\sigma$ is irreducible, so it necessarily follows that $\sigma$ is a
singleton,  say $\sigma=\{s\}$, where $s$ is irreducible.   In turn, the first condition above implies that $s$ lies in no
$F^\theta_r$, whence $Ss=0$.

\definition \label IntroPhiSigmaNonDegMudouNumeros
  \iItemize
  \iItem An element $s$ in $S$ will be called \emph {degenerate} if $s$ is irreducible and $Ss=\{0\}.$
  \iItem A string $\sigma$ will be called \emph {degenerate} if $\sigma=\{s\}$, where $s$ is a degenerate element.
  \iItem The set of all non-degenerate strings will be denoted by $S^\star _\ess $.
  \iItem For every non-degenerate string $\sigma$, we shall denote by $\varphi_\sigma$ the character of $\ehull $ given by
  $$
  \varphi_\sigma(X)= [\sigma\in\varepsilon(X)], \for X\in\ehull .
  $$

Suppose we are given $\varphi_\sigma$ and we want to recover $\sigma$ from $\varphi_\sigma$.  In the special case in which $S$ has right local
units, we have that
  $$
  \varphi_\sigma(E^\theta_s) = 1 \iff
  \sigma\in\varepsilon(E^\theta_s) = E^\star_s \explain {FRFstarR.iii} \iff
  s\in\sigma,
  \equationmark ClueStringRLU
  $$
  so $\sigma$ is recovered as the set    $\{s\in S: \varphi_\sigma(E^\theta_s)=1\}$.  Without assuming right local units, the last part of
\ref {ClueStringRLU} cannot be trusted, but it may be replaced with
  $$
  \phantom {\varphi_\sigma(E^\theta_s) = 1 \iff
  \sigma\in\varepsilon(E^\theta_s) = E^\star_s}
  \cdots \explain {FRFstarR.ii} \iff
  \sigma\cap E^\theta_s\neq\emptyset ,
  \equationmark ClueString
  $$
 so we at least know which $E^\theta_s$ have a nonempty intersection with $\sigma$.

\state Proposition \label PropsInterior
  Given any string $\sigma$, let the \emph {interior} of $\sigma$ be defined  by
  $$
  \interior \sigma := \{s\in S: \exists x\in S,\ sx\in\sigma\}.
  $$
  Then
  $$
  \interior \sigma = \{s\in S: \varphi_\sigma(E^\theta_s)=1\}.
  $$

\medskip Given any character $\varphi$ of $\ehull $, regardless of whether or not it is
of the form $\varphi_\sigma$ as above, we may still consider the set
  $$
  \sigma_\varphi:=\{s\in S: \varphi(E^\theta_s)=1\},
  \equationmark defineSigFi
  $$
  \#  $\sigma_\varphi$; String associated with the character $\varphi$ (when nonempty);
  so that, when $\varphi=\varphi_\sigma$, we get
$\sigma_{\varphi} =\interior \sigma$.

\state Proposition \label StringFromChar
  If $\varphi$ is any character of\/ $\ehull $, and $\sigma_\varphi$ is nonempty, then $\sigma_\varphi$
is a string closed under least common multiples.

Based on \ref {IntroPhiSigmaNonDegMudouNumeros.iv} we may define a map from the set of all non-degenerate strings to
  $\spec $,
  \#  $\textbackslash spec$; The spectrum of $\ehull $;
  the spectrum of $\ehull $, by
  $$
  \Phi: \sigma\in S^\star _\ess \mapsto\varphi_\sigma\in\spec ,
  \equationmark DefinePhi
  $$
  \# $\Phi$; Map sending non-degenerate strings to their associated characters;
  but if we want the dual correspondence suggested by \ref {defineSigFi}, namely
  $$
  \varphi \mapsto \sigma_\varphi,
  \equationmark ProposeSigma
  $$
  to give a well defined map from $\spec $ to $S^\star$, we need to worry about its domain because we have not checked that $\sigma_\varphi$ is always nonempty, and
hence $\sigma_\varphi$ may fail to be a string.  The appropriate domain is evidently given by the set of all characters $\varphi$ such that
$\sigma_\varphi$ is nonempty but, before we formalize this map, it is interesting to introduce a relevant sub{\sla } of $\ehull $.

\state Proposition \label IntroEOne
  The subset of\/ $\ehull $ given by\fn {This should be contrasted with \ref {FormOfHull},  where the general form of an
element of $\ehull $ is $uF^\theta_\Lambda$, where $u$ is in $\tilde S$, rather than $S$.}
  $$
  \eone = \{sF^\theta_\Lambda ,\
  \Lambda\subseteq S \hbox { is finite, and } s\in\Lambda\},
  $$
  \# $\eone $; A sub{\sla } of $\ehull $;
  is an ideal of\/ $\ehull $.
Moreover, for every $X$ in $\ehull $, one has that $X$ lies in $\eone $ if and only if $X\subseteq E^\theta_s$, for some $s$ in $S$.

Whenever $J$ is an ideal in a {\sla} $\E$, there is a standard inclusion
  $$
  \varphi\in \hat J \mapsto \tilde \varphi \in\hat E,
  $$
  where, for every $x$ in $E$,
  one has that  $\tilde \varphi(x) = 1$, if and only if there exists some $y$ in $J$ with $y\leq x$, and $\varphi(y)=1$.
The next result is intended to distinguish the
elements of the  copy of
$\seone $ within $\spec $ given by the above correspondence.

\state Proposition \label DistinguishSeone
  Let $\varphi$ be a character on $\ehull $.  Then the following are equivalent:
  \iItemize
  \iItem $\varphi\in\seone $, 
  \iItem $\varphi(E^\theta_s)=1$, for some $s$ in $S$,
  \iItem $\sigma_\varphi$ is nonempty,  and hence it is a string by \ref {StringFromChar}. 

If $S$ admits right local units, then $\eone =\ehull $, so $\sigma_\varphi$ is a string for every character $\varphi\in\spec $.

The vast majority of non-degenerate strings $\sigma$ lead to a character $\varphi_\sigma$ belonging to $\seone $, but there are exceptions.

\state Proposition \label IrrStrings
  If $\sigma$ is a non-degenerate string in $S^\star _\ess $ then $\varphi_\sigma$ does not belong to $\seone $ if and only if $\sigma=\{s\}$, where $s$ is an
irreducible element of $S$.

By \ref {DistinguishSeone} we have
 that the largest set of characters on which the correspondence described in \ref {ProposeSigma} produces a
bona fide string is precisely $\seone $, so
we may now formaly introduce the map suggested by that correspondence.

\definition \label IntroSigma
  We shall let
  $$
  \Sigma:\seone \mapsto S^\star,
  $$
  \# $\Sigma$; Map sending characters to their associated strings;
  be the map
  given by
  $$
  \Sigma(\varphi) = \sigma_\varphi=\{s\in S: \varphi(E^\theta_s)=1\}, \for \varphi\in\seone .
  $$

For every string $\sigma$,
excluding the exceptional ones discussed in \ref {IrrStrings},  we then have that
  $$
  \Phi(\sigma)=\varphi_\sigma\in\seone ,
  $$
  and
  $$
  \Sigma\big(\Phi(\sigma)\big)=\interior \sigma,
  \equationmark SigmaPhiIsInterior
  $$
  by \ref {PropsInterior}.

\definition A string $\sigma$ in $S$ will be termed \emph{open} if $\sigma=\interior \sigma.$

The nicest situation is for open strings:

\state Proposition \label SigmaPhiOnOpen
  If $\sigma$ is an open string, then
  \iItemize
  \iItem $\sigma$ is non-degenerate,
  \iItem $\Phi(\sigma)\in\seone $, and
  \iItem $\Sigma\big(\Phi(\sigma)\big)=\sigma$.

Given that the composition $\Sigma\circ\Phi$ is so well behaved for open strings, we will now study the reverse composition $\Phi\circ\Sigma$
on a set of characters related to open strings.

\definition \label DefineOpenChar
  A character $\varphi$ in $\spec $ will be called an \emph {open} character if $\sigma_\varphi$ is a (nonempty) open string.

We remark that every open character belongs to $\seone $ by \ref {DistinguishSeone}, although not all characters in $\seone $ are open.

By \ref {SigmaPhiOnOpen} it is clear that $\varphi_\sigma$ is an open character for every open string $\sigma$.

If $S$ admits right local units, we have seen that every string in $S^\star $ is open, and also that
$\sigma_\varphi$ is a string for every character.  Therefore every character in $\spec $ is open.

The composition $\Phi\circ\Sigma$ is not as well behaved as the one discussed in \ref {SigmaPhiOnOpen}, but there is at least
some relationship between a character $\varphi$ and its image under $\Phi\circ\Sigma$, as we shall now see.

\state Proposition \label RioMaior
  Given  any open character $\varphi$, one has that
  $$
  \varphi\leq\Phi\big (\Sigma(\varphi)\big ).
  $$

This leads us to one of our main results.

\state Theorem \label BigNewTightResult
  Let $S$ be a $0$-left-cancellative semigroup admitting least common multiples.  Then, for every open, maximal string
$\sigma$ over $S$, one has that $\varphi_\sigma $ is an ultra-character.

The previous result raises the question as to whether $\sigma_\varphi$ is a
maximal string for every ultra-character $\varphi$, but this is not true in
general.  Consider for example the unital semigroup
  $$
  S=\{1,a,aa,0\},
  $$
  in which $a^3=0$.
  The $\theta$-constructible subsets of $S$ are precisely

  \vbox { \bigskip \begingroup \offinterlineskip
  $$
  \vcenter {\halign {
  \vrule height 14pt depth 8pt\ \ \hfill #\hfill \ \ &\vrule \ \ \hfill #\hfill
\ \ &\vrule \ \ \hfill #\hfill \ \ \vrule \cr \tabrule
  $E^\theta_1 = F^\theta_1 = \{1,a,aa\}$ & & \cr \tabrule
  $\hfill F^\theta_a = \{1,a\}$ & $E^\theta_a = aF^\theta_a = \{a,aa\}$ & \cr
\tabrule
  $\hfill F^\theta_{aa} = \{1\}$ & $\hfill aF^\theta_{aa} = \{a\}$ & $E^\theta
_{aa} = aaF^\theta_{aa} = \{aa\}$ \cr \tabrule
  }}
  $$
  \vskip -5pt\centerline {\eightrm List of $\scriptstyle \theta$-constructible
sets}
  \endgroup
  \bigskip }

  \noindent and there are three strings over $S$, namely

  \bigskip \begingroup \offinterlineskip
  $$
  \vcenter {\halign {
  \vrule height 14pt depth 8pt\ \ \hfill #\hfill \ \ &\vrule \ \ \hfill #\hfill
\ \ &\vrule \ \ \hfill #\hfill \ \ \vrule \cr \tabrule
  $\delta_1=\{1\}$ & $\delta_a=\{1,a\}$ & $\delta_{aa}=\{1,a,aa\}$ \cr
\tabrule
  }}
  $$
  \endgroup \bigskip

  Since the correspondence $s\mapsto\delta_s$ is a bijection from $S'$ to
$S^\star$, we see that $\theta^\star$ is isomorphic to $\theta$, and in
particular the $\theta^\star$-constructible subsets of $S^\star$, listed
below, mirror the $\theta$-constructible ones.

  \vbox {
  \bigskip \begingroup \offinterlineskip
  $$
  \vcenter {\halign {
  \vrule height 14pt depth 8pt\ \ \hfill #\hfill \ \ &\vrule \ \ \hfill #\hfill
\ \ &\vrule \ \ \hfill #\hfill \ \ \vrule \cr \tabrule
  $E^\star_1 = F^\star_1 = \{\delta_1,\delta_a,\delta_{aa}\}$ & & \cr
\tabrule
  $\hfill F^\star_a = \{\delta_1,\delta_a\}$ & $E^\star_a = aF^\star_a =
\{\delta_a,\delta_{aa}\}$ & \cr \tabrule
  $\hfill F^\theta_{aa} = \{\delta_1\}$ & $\hfill aF^\star_{aa} = \{\delta
_a\}$ & $E^\theta_{aa} = aaF^\star_{aa} = \{\delta_{aa}\}$ \cr \tabrule
  }}
  $$
  \endgroup
  \vskip -5pt\centerline {\eightrm List of $\scriptstyle \theta^\star
$-constructible sets}
  }
  \bigskip

Observe that the string $\sigma:=\delta_a=\{1,a\}$ is a proper subset of the
string $\{1,a,aa\}$, and hence $\sigma$ is not maximal.  But yet notice that
$\varphi_\sigma$ is an ultra-character, since $\{\delta_a\}$ is a
  minimal\fn {Whenever $e_0$ is a nonzero minimal element of a {\sla } $E$,
the character $\varphi(e)=[e_0\leq e]$ is an ultra-character.}
  member of $\P (S^\star,\theta^\star)$.
  We thus get an example of

\medskip {\it ``A string $\sigma$ which is not maximal but such that $\varphi
_\sigma$ is an ultra-character.''}  \medskip

\noindent On the other hand, since $\sigma=\sigma_{\varphi_\sigma}$, this
also provides an example of

\medskip {\it ``An ultra-character $\varphi$ such that $\sigma_\varphi$ is
not maximal.''}  \medskip

This suggests the need to single out the strings which give rise to
ultra-characters:

\definition
  We will say that a string $\sigma$ is \emph {quasi-maximal} whenever $\varphi
_\sigma$ is an ultra-character.
  The set of all quasi-maximal strings will be denoted by $S^\propto$.
  \# $S^\propto$; The set of all quasi-maximal strings;

Adopting this terminology, the conclusion of \ref {BigNewTightResult} states
that every open,  maximal string is quasi-maximal.

\state Theorem \label UltraIsPhiSigma
  Let $S$ be a $0$-left-cancellative semigroup admitting least common multiples.  Then,
every open ultra-character on $\ehull $ is of the form $\varphi_\sigma$
for some open, quasi-maximal string $\sigma$.

The importance of quasi-maximal strings evidenced by the last result begs for a
better understanding of such strings.  While we are unable to provide a complete
characterization, we can at least exhibit some further examples beyond the
maximal ones.

To explain what we mean, recalll from \ref {FRFstarR.i} that a string $\sigma$
belongs to some $F^\star_\Lambda$ if and only if $\sigma$ is contained in
$F^\theta_\Lambda$.  It is therefore possible that $\sigma$ is maximal among
all strings contained in $F^\theta_\Lambda$, and still not be a maximal
string.  An example is the string $\{1,a\}$ mentioned above, which is maximal
within $F^\theta_a$, but not maximal in the strict sense of the word.

\state Proposition \label RelativelyMaximalStrings
  Let $\Lambda$ be a nonempty finite subset of $S$ and suppose that $\sigma$
is an open string such that $\sigma\subseteq F^\theta_\Lambda$.  Suppose moreover
that $\sigma$ is maximal among the strings contained in $F^\theta_\Lambda$,
in the sense that for every string $\mu$, one has that
  $$
  \sigma\subseteq\mu\subseteq F^\theta_\Lambda\Imply \sigma=\mu.
  $$
  Then $\varphi_\sigma$ is an ultra-character, and hence $\sigma$ is a
quasi-maximal string.

\def \dualrep {\hat \theta}
\def \ac #1#2{\dualrep _{#1}(#2)}
\def \acinv #1#2{\dualrep _{#1}\inv (#2)}
\def \exist {\exists \kern 0.8pt}

\section Ground characters

In the last section we were able to fruitfully study open characters using strings, culminating with Theorem
\ref {UltraIsPhiSigma}, stating that every open ultra-character is given in terms of a string.
However nothing of interest was said about an ultra-character when it is not open.
The main purpose of this section is thus to obtain some useful information about non-open ultra-characters.  The
main result in this direction is Theorem \ref {NonOpenUltra}, below.

\fix Throughout this section we fix a $0$-left-cancellative semigroup $S$ admitting least common multiples.
For each $s$ in $S$ let
  $$
  \hat F_ s = \{\varphi\in \spec : \varphi(F^\theta_s)=1\}, \and
  \hat E_ s = \{\varphi\in \spec : \varphi(E^\theta_s)=1\},
  $$
  \# $\dualrep $; Dual representation of $S$ on $\textbackslash spec$;
  \# $\hat F_s$; Domain of $\dualrep _s$;
  \# $\hat E_s$; Range of  $\dualrep _s$;
  and for every $\varphi$ in $\hat F_ s$, consider the character   $\ac s\varphi$ given by
  $$
  \ac s\varphi(X) = \varphi\big(\theta_s\inv (E^\theta_s\cap X)\big), \for X\in\ehull .
  $$
Observing that
  $$
  \ac s\varphi(E^\theta_s) = \varphi\big(\theta_s\inv (E^\theta_s)\big) = \varphi(F^\theta_s) = 1,
  \equationmark ActionOnPhiIsChar
  $$
  we see that $\ac s\varphi$ is indeed a (nonzero) character,  and that $\ac s\varphi$ belongs to   $\hat E_ s$.  As a
consequence we get a map
  $$
  \dualrep _ s : \hat F_ s \to \hat E_ s,
  $$
  which is easily seen to be bijective, with inverse given by
  $$
  \acinv s\varphi (X) =
  \varphi\big(\theta_s(F^\theta_s\cap X)\big), \for \varphi\in\hat E_ s, \for X\in\ehull .
  $$
  We may then see each $\dualrep _ s$ as an element of $\I \big(\spec \big)$, and it is not hard to see that
the correspondence
  $$
  \dualrep : s\in S \mapsto \dualrep _ s\in\I (\spec )
  $$
  is a representation of $S$ on $\spec $.

All of this may also be deduced from the fact that any inverse semigroup, such as $\hull $, admits a canonical
representation on the
spectrum of its idempotent {\sla } (see \cite [Section 10]{actions}), and that $\dualrep $ may be obtained as the composition
  $$
  S\arw \theta  \hull \longrightarrow \I \big(\spec \big),
  $$
  where the arrow in the right-hand-side is the canonical representation mentioned above.

\definition We shall refer to $\dualrep $ as the \emph {dual representation} of $S$.

The following technical result gives some useful information regarding the relationship between the dual representation and the representation
  $\rho$ of $\hull $ described in \ref {MapForISG}.

\state Lemma \label CovarPhiLemma
  Given $s$ in $S$, and $\sigma$ in $S^\star _\ess $, one has that
  \iItemize
  \iItem $\varphi_\sigma\in\hat F_ s \IFF \sigma\in F^\star _s$,
  \iItem if the equivalent conditions in (i) are satisfied,  then $\ac s{\varphi_\sigma}=\varphi_{\theta^\star _s(\sigma)}$,
  \iItem $\varphi_\sigma\in\hat E_ s \IFF \sigma\in E^\star _s$,
  \iItem if the equivalent conditions in (iii) are satisfied,  then $\acinv s{\varphi_\sigma}=\varphi_{\theta^{\star -1}_s(\sigma)}$.

Considering the representation $\theta^\star $ of $S$ on $S^\star $, observe that $S^\star _\ess $ is an invariant
  subset of $S^\star $,  and it is easy to see that it is also invariant under the representation
  $\rho$ of $\hull $ described in \ref {MapForISG}.  Together with the dual representation of $\hull $ on $\spec $ mentioned
above, we thus have two natural representations of $\hull $, which are closeely related, as the following immediate
consequence of the above result asserts:


\state Proposition \label CovarPhi
  The mapping
  $$
  \Phi:S^\star _\ess \to\spec
  $$
  of \ref {DefinePhi} is covariant relative to the natural representations of $\hull $ referred to above.

The fact that the correspondence between strings and characters
  is not a perfect one
  (see e.g.~\ref {SigmaPhiIsInterior} and \ref {RioMaior})
  is partly responsible for the fact that expressing the  covariance properties of the map $\Sigma$ of  \ref {IntroSigma} cannot be done
in the same straightforward way as we did for $\Phi$ in \ref {CovarPhi}.  Nevertheless, there are some things we may say in
this respect.

Let us first treat the question of covariance regarding
$\acinv s\varphi$.  Of course, for this to be a well defined character we need $\varphi$ to be in $\hat E_ s$, meaning that $\varphi(E^\theta_s)=1$, which is also equivalent to saying that $s\in\sigma_\varphi$.   In particular characters with empty strings are immediately
ruled out.

\state Lemma \label BackInvarString
  For every $s$ in $S$, and every character $\varphi$ in $\hat E_ s$, one has that
  $$
  \sigma_{\acinv s\varphi} = \{p\in S: sp\in\sigma_\varphi\}.
  $$

The set appearing in the right hand side of the equation displayed in  \ref {BackInvarString} is precisely the same set
mentioned in definition \ref {MappingStrings.ii}
of $s\inv *\sigma_\varphi$, except that this notation is reserved for the situation in which the intersection of $\sigma$ with $sS$ is nonempty, which
precisely means that $s\in\interior \sigma_\varphi$.

\state Proposition \label DualBackOnStrings
  Pick $s$ in $S$ and let $\varphi$ be any  character in $\hat E_ s$.  Then $s\in\sigma_\varphi$, and moreover
  \iItemize
  \iItem if $s$ is in $\interior \sigma_\varphi$, then $\sigma_\varphi\in E^\star _s$, and $\sigma_{\acinv s\varphi} = \theta^{\star -1}_s(\sigma_\varphi)$,
  \iItem if $s$ is not in $\interior \sigma_\varphi$, then $\sigma_{\acinv s\varphi} = \emptyset $.

Regarding the behavior of strings associated to characters of the form $\ac s\varphi$, we have:

\state Lemma \label BirthOfString
  For every $s$ in $S$, and every character $\varphi$ in $\hat F_ s$, one has that $\ac s\varphi$ belongs to $\seone $ (and hence
\ref {DistinguishSeone} implies that
$\sigma_{\ac s\varphi}$ is a string), and moreover
  \iItemize
  \iItem if $\sigma_\varphi$ is nonempty, then $\sigma_\varphi\in F^\star _s$, and      $\sigma_{\ac s\varphi} =   \theta^\star _s(\sigma_\varphi)$,
  \iItem if $\sigma_\varphi$ is empty, then $\sigma_{\ac s\varphi} =     \delta_s$.

We may interpret the above result, and more specifically the identity
  $$
  \sigma_{\ac s\varphi} =   \theta^\star _s(\sigma_\varphi),
  $$
  as saying that the correspondence $\varphi\mapsto\sigma_\varphi$ is covariant with respect to the actions
$\dualrep $ and $\theta^\star $, on $\spec $ and $S^\star $, respectively, except that the term ``$\sigma_\varphi$" appearing  is the
right-hand-side above is not a well defined string since it may be empty, even though the left-hand-side  is always well defined.
In the problematic case of an empty string, \ref {BirthOfString.ii} then gives the undefined right-hand-side the default value of $\delta_s$.

	%
	%
	%


\definition A character $\varphi$ in $\spec $ will be called a \emph {ground} character if $\sigma_\varphi$ is empty.

By \ref {DistinguishSeone},  the ground characters  are precisely the members of $\spec \setminus \seone $.

Besides the ground characters, a character $\varphi$ may fail to be open because $\sigma_\varphi$, while being a bona fide string,  is
not an open string.  In this case we have
  that $\sigma_\varphi=\delta_s$, for some $s$ in $S$ such that $s\notin sS$.

\state Proposition \label PrepOrbitGroundChar
  Let $\varphi$ be a character such that $\sigma_\varphi=\delta_s$, where $s$ is such that $s\notin sS$.  Then $\varphi\in\hat E_ s$, and $\acinv s\varphi$ is a
ground character.

We may now give  a precise characterization of non-open characters in terms of the ground characters:

\state Proposition \label NonOpenChars
  Denote by \
  $\sop $
  \# $\sop $; Set of all open characters on $\ehull $;
  the set of all open characters on $\ehull $.  Then
  $$
  \spec \setminus \sop = \big \{\ac u\varphi: u\in\tS , \ \varphi \hbox { is a ground character in } \hat F_ u\big \}.
  $$
  Moreover for each $\psi$ in the above set, there is a unique pair $(u, \varphi)$, with $u$ in $\tS $, and $\varphi$ a ground character,
such that $\psi=\ac u\varphi$.

We may now combine several of our earlier results to give a description of all ultra-characters on $\ehull $.

\state Theorem \label NonOpenUltra
  Let $S$ be a $0$-left-cancellative semigroup admitting least common multiples.  Denote by
  $\sinf $
  \# $\sinf $; Set of all ultra-characters on $\ehull $;
  the set of all ultra-characters on $\ehull $, and by
  $$
  \sinfop = \sop \cap \sinf ,
  $$
  \# $\sinfop $; Set of all open ultra-characters  on $\ehull $;
  namely  the subset formed by all open ultra-characters.  Then
  \iItemize
  \iItem $\sinfop \phantom {\setminus \sinfop \ } = \big \{\varphi_\sigma: \sigma \hbox { is an open, quasi-maximal string in }S\big \}$,  and
  \iItem $\sinf \setminus \sinfop = \big \{\ac u\varphi: u\in\tS , \ \varphi \hbox { is a ground, ultra-character in } \hat F_ u\big \}$.

The upshot is that in order to understand all ultra-characters on $\ehull $,  it only remains to describe the ground ultra-characters.

\section Subshift semigroups

\label SubshiftSgrSection
By a subshift on a finite alphabet  $\Sigma$
  one means a subset $\X \subseteq\Sigma^{\bf N}$,
  which is closed relative to the product topology, and invariant under the left shift map
  $$
  x_1x_2x_3\ldots \ \mapsto \ x_2x_3x_4\ldots
  $$

\fix Throughout this chapter we will let $\X $ be a fixed subshift. \medskip

The language of $\X $ is the set $L$ formed by all finite words appearing as a block in some infinite word belonging to
$\X $.  We will not allow the empty word in $L$, as sometimes done in connection with subshifts, so all of our words have
strictly positive length.

In the present section we will be concerned
with the semigroup
  $$
  \SX =L\cup\{0\},
  \equationmark IntroSX
  $$
  equipped with the multiplication operation given by
  $$
  \def \quad { }
  \mu \cdot \nu =\left \{\matrix {
  \mu \nu , & \hbox { if $\mu ,\nu \neq 0$, and $\mu \nu \in L$},
\cr \pilar {12pt}
  0, & \hbox { otherwise, }\hfill }\right .
  $$
  where $\mu\nu$ stands for the concatenation of $\mu$ and $\nu$.

Given $\mu$ and $\nu$ in $\SX $, with $\nu$ nonzero, notice that $\mu\|\nu$ if and only if $\mu$ is a prefix of $\nu$.  Given that
divisibility is also well defined for the unitized semigroup $\tilde \SX $, and that $1$ divides any $\mu\in\SX $, we
will also say that $1$ is a prefix of $\mu$.

Some of the special properties of $\SX $ of easy verification are listed below:

\state Proposition  \label EasyProp
  \iItemize
  \iItem $\SX $ is $0$-left-cancellative and $0$-right-cancellative,
  \iItem $\SX $ admits least common multiples,
  \iItem $\SX $ has no idempotent elements other than $0$.

A further special property of $\SX $ is a very strong uniqueness of the normal form for elements in $\hullX $:

\state Proposition \label SubshiftUniqueness
  For $i=1,2$, let $\Lambda_i$ be a finite subset of $\tSX $ intersecting $\SX $ non-trivially, and let $u_i,v_i\in\Lambda_i$ be such that
  $$
  \theta_{u_1}f_{\Lambda_1} \theta_{v_1}\inv = \theta_{u_2}f_{\Lambda_2} \theta_{v_2}\inv \neq0.
  $$
  Then $u_1=u_2$,  $v_1=v_2$, and $F_{\Lambda_1} = F_{\Lambda_2}$.

Given the importance of strings, let us give an explicit description of these in the present context.

\state Proposition \label DescribeStringsAsWords
  Given a (finite or infinite) word
  $$
  \omega = \omega_1\omega_2\omega_3\ldots,
  $$
  on the alphabet $\Sigma$, assume $\omega$  to be \emph {admissible} (meaning that $\omega$ belongs to $L$, if finite, or to $\X $, if
infinite) and
consider the set $\delta_\omega$ formed by all prefixes of $\omega$ having positive length, namely
  $$
  \delta_\omega=\{\omega_1,\ \omega_1\omega_2,\  \omega_1\omega_2\omega_3,\ \ldots\ \}.
  $$
  Then:
  \iItemize
  \iItem $\delta_\omega$ is a string,
  \iItem $\delta_\omega$ is an open string if and only if $\omega$ is an infinite word,
  \iItem $\delta_\omega$ is a maximal string if and only if $\omega$ is an infinite word,
  \iItem for any string $\sigma$ in $\SX $,  there exists a unique admissible word  $\omega$ such that $\sigma=\delta_\omega$.

Strings consist in one of our best instruments to provide characters on $\ehullX $. Now that we have a concrete
description of strings in terms of admissible words, let us give an equally concrete description of the characters
induced by strings.

\state Proposition \label DescribeCharsAsWords
  Let $\omega$ be a given  infinite admissible word, and let $X$ be any $\theta$-constructible set, written in normal form,
namely $X=uF^\theta_\Lambda$,
where $\Lambda$ is a finite subset of $\tSX $, intersecting $\SX $ nontrivially, and $u\in\Lambda$.  Regarding the string $\delta_\omega$, and
the associated character $\varphi_{\delta_\omega}$, the following are equivalent:
  \iItemize
  \iItem $\varphi_{\delta_\omega}(X)=1$,
  \iItem $u$ is a prefix of $\omega$, and
upon writing $\omega=u\eta$, for some infinite word $\eta$, one has that $t\eta$ is admissible (i.e.~lies in $\X $), for every $t$ in $\Lambda$.

In case the set $X$ of the above result coincides with $F^\theta_\mu$,  for some $\mu$ in $L$,  we get the following
simplification:

\state Proposition \label DescribeCharsAsWordsOnFollower
  Let $\omega$ be an infinite admissible word, and let $\mu\in L$.  Then
  $$
  \varphi_{\delta_\omega}(F^\theta_\mu)=1 \IFF \mu\omega\in\X .
  $$

Let us now study a situation in which infinite words provide all ultra-characters.

\state Proposition \label GroundForSubshifts
  The following are equivalent:
  \iItemize
  \iItem for every finite subset $\Lambda\subseteq\tSX $, such that $\Lambda\cap\SX \neq\emptyset $, one has that $F^\theta_\Lambda$ is either empty or infinite,
  \iItem every nonempty constructible set is infinite,
  \iItem $\{E^\theta_a: a\in\Sigma\}$ is a cover for $\ehullX $,
  \iItem $\ehullX $ admits no ground ultra-characters,
  \iItem for every ultra-character $\varphi$ on $\ehullX $, there exists an infinite admissible word $\omega$, such that $\varphi=\varphi_{\delta_\omega}$.

Let us conclude this section with an example to show that nonempty finite constructible sets may indeed exist and hence
the equivalent conditions of  \ref {GroundForSubshifts} do not
always hold.  Consider the alphabet $\Sigma =\{a, b, c\}$ and let $\X $ be the subshift on $\Sigma $ consisting of all infinite words $\omega$
such that, in any block of $\omega$ of length three, there are no repeated letters.  Alternatively, a set of forbidden words
defining $\X $ is the set of all words of length three with some repetition.

It is then easy to see that the language $L$ of $\X $ is formed by all finite words on $\Sigma $  with the same restriction on
blocks of length three described above.

Notice that $c\in F_{\{a, b\}}$, because both $ac$ and $bc$ are in $L$.  However there is no element in $F_{\{a, b\}}$
other than $c$, because it is evident that neither $a$ nor $b$ lie in $F_{\{a,b\}}$, and for any $x$ in $\Sigma $,
either $acx$ or $bcx$ will involve a repetition. So \emph {voil\`a} the finite constructible set:
  $$
  F_{\{a,b\}}=\{c\}.
  \equationmark FiniteConstructibleSet
  $$

\def\GS {{\cal G}_\S}

\section C*-algebras associated to subshifts

In this final section we will briefly discuss applications of the theory so far described to various C*-algebras
associated to subshifts that have been studied starting with Matsumoto's original work \cite{MatsuOri}.

Given a subshift $\X$, as in the previous section, we will consider the $0$-left-cancellative semigroup $\SX$, as well
as its inverse hull $\hullX$.  We may then consider several general constructions of C*-algebras from inverse
semigroups, and our goal is to argue that many of these, once applied to $\hullX$, produce all of the C*-algberas
studied in the literature in connection with subshifts.

The constructions we have in mind share a common pattern in the following sense.  Given an inverse semigroup $\S$ with
zero, consider the standard action of $\S$ on $\widehat E (\S)$, namely the dual of the idempotent {\sla} of $\S$.  We
may then build the groupoid $\GS$ formed by all germs for this action.  This groupoid is sometimes  referred to via the
suggestive notation
  $$
  \S \ltimes \widehat E (\S).
  $$
  The C*-algebra of $\GS$ is well known to be a
  quotient\fn{Modulo the relation that identifies the zero of $\S$ with the zero of the corresponding C*-algebra.}
  of Paterson's universal C*-algebra for $\S$.

Whenever $Y$ is a closed invariant subset of $\widehat E (\S)$, we may restrict the action of $\S$ to $Y$,  and consider the
corresponding groupoid of germs
  $\S \ltimes Y$,
  which may also be seen as the reduction of $\GS$ to $Y$.

Our point is that several C*-algebras studied in the
literature in connection to the subshift $\X$ are actually groupoid C*-algebras of the form
  $$
  C^*\big(\hullX \ltimes Y\big),
  $$
  where $Y$ is a closed subset of $\specX$, invariant under the standard action of $\hullX$.

In order to describe the first relevant alternative for $Y$, let $S$ be any
$0$-left-cancellative semigroup and consider the representation $\iota$ of $\ehull$ on $\P(S')$ given by the inclusion of the former in
the latter.  We will say that a character $\varphi$ of $\ehull$  is essentially tight (relative to the above representation $\iota$)
provided one has that
  $$
  \varphi(X)=\bigvee _{i=1}^n\varphi(Y_i),
  $$
  whenever $X, Y_1,\ldots,Y_n$ are in $\ehull $, and the symmetric
difference
  $$
  X\mathop {\Delta}\big(\bigcup _{i=1}^nY_i\big)
  $$
  is finite.  The set of all essentially tight characters of $\ehull$ will be denoted by $\sess$, and it may be shown
that $\sess$ is a closed invariant subset of $\spec$.

The second relevant alternative for $Y$ is based on the set $\smax$ defined by
  $$
  \smax = \{\varphi_\sigma: \sigma \hbox { is a maximal string}\}.
  $$
  In general   $\smax$ is not invariant under the action of $\hull$, but when $\S=\SX$ for some subshift $\X$,
invariance is guaranteed.  One may then take $Y$ to be the closure of
  \def\varspecX#1{\gothEhat_\sixrmbox{#1} (\SX)}%
  \def\G#1{{\cal G}_\sixrmbox{#1}}
  $\varspecX{max}$.

The tight spectrum  $\gothEhat_\tight (\SX)$ is a further alternative,  and in some sense it is the most natural one given
that many C*-algebras associated to inverse semigroups turn out to be the groupoid C*-algebra for the reduction of
Paterson's universal groupoid to the tight spectrum of the idempotent semilattice.

These subsets are related to each other as follows
  $$
  \def \quad {\kern 8pt}
  \matrix {
  \overline{\varspecX{max}} & \subseteq & \varspecX{tight} \cr
  \kern 15pt \rotatebox {-90}{$\subseteq$}\hfill \cr \pilar {16pt}
  \varspecX{ess}}
  \equationmark RelationSets
  $$

As already observed, all of these are subsets of $\ehullX$ which are closed and
invariant under the action of $\hullX$, so each gives rise to a reduced subgroupoid, which we will correspondingly denote by $\G{max}$,
$\G{tight}$ and $\G{ess}$.

\state Theorem
  Given any subshift $\X$, one has that
  \iItemize
  \iItem $C^*(\G{ess})$ is isomorphic to Matsumoto's C*-algebra introduced in \cite{MatsuOri}.
  \iItem $C^*(\G{max})$ is isomorphic to Carlsen-Matsumoto's C*-algebra introduced in \cite[Definition 2.1]{MatsuCarl}.

There are many situations in which the inclusions in \ref{RelationSets} reduce to equality, but examples may be given to
show these are, in general, proper inclusions.  The fact that $\overline{\varspecX{max}}$ and $\varspecX{ess}$ may
differ is related to the fact that Matsumoto's C*-algebra may be non-isomorphic to the Carlsen-Matsumoto one, but it may
be shown that under condition ($*$)\fn{The reader should be warned that the description of condition ($*$) in
\cite{MatsuCarl} is incorrect and must be amended by requiring that the sequence $\{\\mu_i\}_i$, mentioned there, have an
infinite range.} of \cite {MatsuCarl}, one has that $\overline{\varspecX{max}}=\varspecX{ess}$, whence isomorphism
holds.

\color{black}

\references

\Article MatsuCarl
  T. M. Carlsen and K. Matsumoto;
  Some remarks on the C*-algebras associated with subshifts;
  Math. Scand., 95 (2004), 145-160

\Article Cherubini
  A. Cherubini and M. Petrich;
  The Inverse Hull of Right Cancellative Semigroups;
  J. Algebra, 111 (1987), 74-113

\Bibitem CP
A. H. Clifford and G. B. Preston;
The algebraic theory of semigroups. Vol. I;
Mathematical Surveys, No. 7.
 American Mathematical Society, Providence, R.I., 1961

\Article CoOne
  L. A. Coburn;
  The C*-algebra generated by an isometry I;
  Bull. Amer. Math. Soc., 73 (1967), 722-726

\Article CoTwo
  L. A. Coburn;
  The C*-algebra generated by an isometry II;
  Trans. Amer. Math. Soc., 137 (1969), 211-217

\Article DokuchaExel
  M. Dokuchaev and R. Exel;
  Partial actions and subshifts;
  J. Funct. Analysis, 272 (2017), 5038-5106

\Article actions
  R. Exel;
  Inverse semigroups and combinatorial C*-algebras;
  Bull. Braz. Math. Soc. (N.S.), 39 (2008), 191-313

\Article infinoa
  R. Exel and M. Laca;
  Cuntz-Krieger algebras for infinite matrices;
  J. reine angew. Math., 512 (1999), 119-172

\Article ExelPardo
  R. Exel and E. Pardo;
  The tight groupoid of an inverse semigroup;
  Semigroup Forum, 92 (2016), 274-303

\Article Li
  X. Li;
  Semigroup C*-algebras and amenability of semigroups;
  J. Funct. Anal., 262 (2012), 4302-4340

\Article MatsuOri
  K. Matsumoto;
  On C*-algebras associated with subshifts;
  Internat. J. Math., 8 (1997), 357-374

\Article MH
M. Morse and G. A. Hedlund;
Unending chess, symbolic dynamics and a problem in semigroups;
Duke Math. J., 11 (1944), 1-7.

\Article Munn
  W. D. Munn;
  Brandt congruences on inverse semigroups;
  Proc. London Math. Soc., {\rm (3)} 14 (1964), 154-164.

\Article MurOne
  G. J. Murphy;
  Ordered groups and Toeplitz algebras;
  J. Operator Theory,   18 (1987), 303-326

\Article MurTwo
  G. J. Murphy;
  Ordered groups and crossed products of C*-algebras;
  Pacific J. Math., 2 (1991), 319-349

\Article MurThree
  G. J. Murphy;
  Crossed products of C*-algebras by semigroups of automorphisms;
  Proc.  London Math. Soc., 3 (1994), 423-448

\Article Nica
  A. Nica;
  C*-algebras generated by isometries and Wiener-Hopf operators;
  J. Operator Theory, 27 (1992), 17-52

\Bibitem Paterson
  A. L. T. Paterson;
  Groupoids, inverse semigroups, and their operator algebras;
  Birkh\umlaut auser, 1999

\endgroup

\bye